 \documentclass[11pt]{article}
 \usepackage[top=1in,bottom=1in,left=1.5in,right=1.5in]{geometry}
  \usepackage{amsmath,amssymb}
  \usepackage{latexsym}
  \usepackage[dvips]{pstcol} 
  \usepackage[dvips]{graphicx}
  \newcommand{\twocases}[5]{#1 =
               \left\{
                \begin{alignedat}{2}
                 &#2 \quad & &\text{if } #3\\
                 &#4 \quad & &\text{if } #5
                \end{alignedat}
               \right.}
  \newcommand{\stack}[2]{\genfrac{}{}{0pt}{}{#1}{#2}}
  \newcommand{\C}{\mathbb{C}}
  \newcommand{\N}{\mathbb{N}}
  \newcommand{\R}{\mathbb{R}}
  \newcommand{\Z}{\mathbb{Z}}
  \renewcommand{\i}{\mathbf{i}}
  \renewcommand{\j}{\mathbf{j}}
  \newcommand{\e}{\mathbf{e}}
  
  \newcommand{\m}{\mathbf{m}}
  \newcommand{\n}{\mathbf{n}}
  \newcommand{\p}{\mathbf{p}}

  \renewcommand{\v}{\mathbf{v}}
  \newcommand{\w}{\mathbf{w}}
  \newcommand{\x}{\mathbf{x}}
  \newcommand{\y}{\mathbf{y}}
  \newcommand{\z}{\mathbf{z}}
  \newcommand{\1}{\mathbf{1}}
  \newcommand{\bM}{\mathbf{M}}

  \newcommand{\cE}{\mathcal{E}}
  \newcommand{\cG}{\mathcal{G}}
  
  \newcommand{\cO}{\mathcal{O}}
  \newcommand{\cP}{\mathcal{P}}

  \newcommand{\cX}{\mathcal{X}}
  \newcommand{\lan}{\langle}
  \newcommand{\ran}{\rangle}
  \newcommand{\an}[1]{\lan#1\ran}
  \newcommand{\hs}{\hspace*{\parindent}}
  \newcommand{\proof}{\hs \textbf{Proof.\ }}
  \newcommand{\tr}{\mathop{\mathrm{tr}}\nolimits}
  \newcommand{\Aut}{\mathop{\mathrm{Aut}}\nolimits}
  \newcommand{\argmax}{\mathop{\mathrm{arg\,max}}}
  \newcommand{\trans}{^\top}
  \newcommand{\qed}{\hspace*{\fill} $\Box$\\}
  \newcommand{\per}{\mathop{\mathrm{perm}}\nolimits}
  \newcommand{\perio}{\mathrm{per}}
  \newcommand{\pr}{\mathrm{pr}}
  \newcommand{\pers}{\per_s}
  \newcommand{\spec}{\mathrm{spec}}
  \newtheorem{theo}{\bfseries \hs Theorem}[section]
  \newtheorem{defn}[theo]{\bfseries \hs Definition}
  \newtheorem{prop}[theo]{\bfseries \hs Proposition}
  \newtheorem{lemma}[theo]{\bfseries \hs Lemma}
  \newtheorem{corol}[theo]{\bfseries \hs Corollary}

  \numberwithin{equation}{section} 

\begin{document}

 \title{Theory of Computation of Multidimensional Entropy with an Application to the Monomer-Dimer Problem}

 \author{
 Shmuel Friedland\thanks{Additional affiliation: New Directions Visiting Professor, Institute of Mathematics
 and its Applications, University of Minnesota, Minnesota, MN 55455-0436.}\\
 \texttt{friedlan@uic.edu}
                  \and
 Uri N. Peled\thanks{This author thanks the Caesarea Edmond Benjamin de Rothschild Foundation
 Institute for Interdisciplinary Applications of Computer Science
 at the University of Haifa, Israel, for partial support.}\\
 \texttt{uripeled@uic.edu}
 }
 \date{Department of Mathematics, Statistics, and Computer Science,\\
       University of Illinois at Chicago\\
       Chicago, Illinois 60607-7045, USA\\
       February 1, 2004}

 \maketitle

 \begin{abstract}
 We outline the most recent theory for the computation of the
 exponential growth rate of the number of configurations on a
 multi-dimensional grid. As an application we compute the
 monomer-dimer constant for the $2$-dimensional grid to $8$
 decimal digits, agreeing with the heuristic computations of
 Baxter, and for the $3$-dimensional grid with an error smaller than
 $1.35\%$.\\[\baselineskip]
 \end{abstract}

 \section{Introduction}

 The exponential growth rate $h$ (with respect to the natural
 logarithm) of the number of configurations on a multi-dimensional
 grid arises in the theory of various phenomena \cite{Pau,FoR}. In
 physics $e^h$ is viewed as the entropy (per atom) of the
 corresponding ``hard model''; in mathematics $h$ is called the
 \emph{topological entropy} \cite{Fr1}; and in information theory
 $h$ (with respect to $\log_2$) is called the
 \emph{multi-dimensional capacity} \cite{WeB}.  The
 $1$-dimensional case is easy, namely $e^h$ is equal to the
 spectral radius $\rho(A)$ of a certain matrix $A$ called the
 ``transfer matrix''. There are very few $2$-dimensional models
 where the value of $h$ is known in closed form
 \cite{Fis,Kas,L,Lie,Bax}.  In all other cases there are estimates
 based on: (a) asymptotic expansions, e.g., \cite{Na2,Bax1,Gau};
 (b) Monte-Carlo methods, e.g., \cite{Ha3,BeS}; (c) bounds, e.g.,
 \cite{Ha2,Ciu,Lun,CaW,FoJ,NaZ}.  In what follows we give a
 complete up-to-date theory of the computation of $h$ by using
 lower and upper bounds.  It refines the techniques described in
 \cite{Fr2} by using an automorphism subgroup of a given graph. A
 fundamental problem in lattice statistics is the monomer-dimer
 problem (see \cite{KeRaSi}). As a demonstration of our
 techniques, we compute the topological entropy of the
 monomer-dimer covers of the $2$-dimensional grid $h_2 =
 .66279897$ (which agrees with the heuristic estimation $e^{h_2} =
 1.940215351$ due to Baxter \cite{Bax1}) and of the $3$-dimensional
 grid $.7653 \leq h_3 \leq .7862$. These numerical results are
 much better than previously known ones.

 Consider the grid $\Z^d$ in $d$-dimensional space $\R^d$.  At
 each point of the grid we place an element of a set of $n$ kinds
 of colors (atoms) denoted by $\an{n} := \{1,\ldots ,n\}$. Certain
 restrictions may be imposed on the colorings. For example, the
 restrictions of the \emph{hard model} are specified by a directed
 $d$-graph $\Gamma:=(\Gamma_1,\ldots ,\Gamma_d)$ called a
 \emph{nearest neighbor digraph}, with $\Gamma_k \subseteq \an{n}
 \times \an{n}$, in the sense that two atoms of kinds $p$ and $q$
 are allowed to occupy respectively the adjacent grid points
 $\i=(i_1,\ldots ,i_d)$ and $\i + \e_k$ (where $\e_k :=
 (\delta_{1k},\ldots ,\delta_{dk})$) only if $(p,q) \in \Gamma_k$.
 We call such a placement a \emph{$\Gamma$-configuration} or
 \emph{$\Gamma$-cover}. This general model is anisotropic, since
 the $\Gamma_k$ can be distinct. A digraph $\Gamma_k$ is called
 \emph{symmetric} when $(p,q)\in \Gamma_k \Leftrightarrow (q,p)\in
 \Gamma_k$. We call $\Gamma$ a \emph{symmetric isotropic nearest
 neighbor digraph} when $\Gamma_1 = \cdots = \Gamma_d = \Delta$, and
 $\Delta$ is symmetric. Let $\m = (m_1,\ldots ,m_d) \in \N^d$,
 where $\N := \{1,2,\ldots\}$, and consider the box $\an{\m} :=
 \an{m_1} \times \cdots \times \an{m_d}$ of dimensions $m_1,
 \ldots ,m_d$.  Let $W(\m)$ be the set of all
 $\Gamma$-configurations of $|\m|_{\pr}:= m_1 \cdots m_d$ atoms in
 the box $\lan \m\ran$. It is easy to show that the sequence $\log
 \#W(\m)_{\m \in \N^d}$ is \emph{subadditive} in each coordinate,
 i.e., $\log \#W(\m + p \e_k) \leq \log \#W(\m) + \log \#W(\m + (p
 - m_k) \e_k)$ for all $\m \in \N^d$, $p \in \N$ and $k \in
 \an{d}$. From this it can be shown that the following limit
 exists and is non-negative or equal to $-\infty$ (we use $\m \to
 \infty$ as an abbreviation of $m_1,\ldots ,m_d \to \infty$):
 \begin{equation}\label{hGamma}
  h = h(\Gamma) := \lim_{\m \to \infty} \frac{\log \#W(\m)}{|\m|_{\pr}},
 \end{equation}
 and  each $\m \in \N^d$ satisfies
 \begin{equation}\label{ub1}
 h \leq \frac{\log \#W(\m)} {|\m|_{\pr}}.
 \end{equation}
 The limit $h(\Gamma)$ is the exponential growth rate of $\#W(\m)$
 per atom, also called \emph{entropy} or \emph{Shannon capacity}.
 It follows from K\"{o}nig's Infinity Lemma that $h = \log
 0=-\infty$ if and only if there are no $\Gamma$-covers of $\Z^d$.
 The case $d=1$ is well understood: $h=\log\rho(A)$, where $A$ is
 the incidence matrix for the digraph $\Gamma_1$; there exist
 $\Gamma$-covers if and only if $\Gamma_1$ has a directed cycle,
 and in that case $h$ is also the exponential growth rate per atom
 of the number of periodic $\Gamma$-covers of $\Z$ \cite{Fr1}. A
 periodic $\Gamma$-cover of $\Z^d$ with period $\m$ (i.e., a
 $\Gamma$-cover $\phi = (\phi_{\i})_{\i \in \Z^d}$ satisfying
 $\phi_{\i + m_k\e_k} = \phi_{\i}$ for all $\i \in \Z^d$ and $k
 \in \an{d}$) is equivalent to a $\Gamma$-cover of the torus
 $T(\m):=(\Z/m_1\Z) \times \cdots \times (\Z/m_d\Z)$.  For $d \geq
 2$, the question whether there exist $\Gamma$-covers is
 undecidable and $h$ is not computable in general \cite{Ber,HKC} (we
 say that a quantity $Q$ is \emph{computable} when given $\epsilon
 > 0$, we can find in a finite number of steps, depending on
 $\epsilon$, a rational number $r$ satisfying $|Q-r| < \epsilon$).
 Equivalently, there is a $d$-digraph $\Gamma$ for which there are
 $\Gamma$-covers of $\Z^d$ but none is periodic. Hence there are
 no nontrivial lower bounds for $h$ in this case. A fundamental
 result in \cite{Fr1} asserts that if at least $d-1$ digraphs out of
 $\Gamma_1, \ldots ,\Gamma_d$ are symmetric, then the exponential
 growth rate per atom of the number of periodic configurations is
 equal to $h$ and $h$ is computable, i.e., we have lower bounds on
 $h$ that converge to $h$.  For $d=2,3$ this will also follow from
 our results in Section 3. In particular these results hold for a
 symmetric isotropic nearest neighbor digraph.

 We mention briefly the topological entropy. Let $W_{\text{top}}(\m)$ be the
 set of all distinct restrictions of $\Gamma$-covers of $\Z^d$ to
 the box $\an{\m}$. $\log \# W_{\text{top}}(\m)$ is also subadditive, and the
 \emph{topological entropy} of $\Gamma$ is defined by
 \[h_{\text{top}}(\Gamma) := \lim_{\m \to \infty} \frac{\log \#W_{\text{top}}(\m)}{|\m|_{\pr}}.\]
 Since $W_{\text{top}}(\m) \subseteq W(\m)$, we have $h_{\text{top}}(\Gamma)
 \leq h(\Gamma)$; a result in \cite{Fr1} asserts that equality
 holds.

 We now elaborate our results. Fix $\m':=(m_1,\ldots ,m_{d-1}) \in
 \N^{d-1}$ and let $\Gamma':=(\Gamma_1,\ldots ,\Gamma_{d-1})$. Let
 $\Omega_d(\m')$ be the transfer digraph between $\Gamma'$-covers
 of $\an{\m'}$ with respect to $\Gamma_d$. That is, the vertex set
 of $\Omega_d(\m')$ is the set of $\Gamma'$-covers of $\an{\m'}$,
 and vertices $u,v$ satisfy $(u,v) \in \Omega_d(\m')$ if and only if
 $[u,v] \in W(\m',2)$, where $[u,v]$ is the configuration
 consisting of $u,v$ occupying the levels $x_d=1,2$ of
 $\an{(\m',2)}$, respectively. We show that $h \leq
 \frac{\log\rho(\Omega_d(\m'))}{|\m'|_{\pr}}$ (by definition,
 the spectral radius of a digraph is the spectral radius of its
 incidence matrix). When $\Gamma_1,\ldots , \Gamma_{d-1}$ are
 symmetric, this upper bound can be improved as follows. Let
 $\Theta_d(\m')$ be the induced subdigraph of $\Omega_d(\m')$ whose
 vertices are the periodic $\Gamma'$-covers of $\an{\m'}$ with
 period $\m'$.  Then we show \cite{Fr2}
 \begin{equation}\label{entdef}
 h(\Gamma) \leq \frac{\log\rho(\Theta_d(\m'))}{|\m'|_{\pr}},
 \;\; m_1,\ldots ,m_{d-1} \text{ even, }
 \Gamma_1,\ldots,\Gamma_{d-1} \text{ symmetric}.
 \end{equation}
 Furthermore, for $\Gamma_1,\ldots,\Gamma_{d-1}$ symmetric, we
 give various lower bounds on $h$ in terms of
 $\log\rho(\Theta_d(\m'))$ for various values of $\m'$.  For
 example, for $d=2$ we show that when $\Gamma_1$ is symmetric, $h
 \geq \frac{\log\rho(\Theta_2(p+2q))-\log\rho(\Theta_2(2q))} {p}$
 for all $p \in \N$ and $q \in \Z_+ := \N \cup \{0\}$. See
 \cite{Fr2} for slightly different lower bounds on $h$, which do
 not use periodicity.

 All of these upper and lower bound converge to the true entropy
 when $\m' \to \infty$.

 We can enhance the efficiency of computing the spectral radius
 $\rho(\Lambda)$ of a digraph $\Lambda \subseteq N \times N$ as
 follows. To compute $\rho(\Lambda)$ one needs to compute the
 spectral radius of its $0$-$1$ $N \times N$ incidence matrix $A$.
 Suppose that $\cG \subseteq S_N$ is an automorphism subgroup of
 $\Lambda$.  Let $\cO=\an{N}/\cG$ be the orbit space under the
 action of $\cG$ and set $M=\#\cO$. Let $\Lambda' \subseteq
 \cO \times \cO$ be the multidigraph induced by $\Lambda$ and $\cG$.
 That is, for $\alpha,\beta \in \cO$, the multiplicity of the edge
 $(\alpha,\beta)$ of $\Lambda'$ is $\widehat{a}_{\alpha,\beta} =
 \sum_{j \in \beta} a_{i,j}$ for any $i \in \alpha$. We show that
 $\rho(\Lambda)$ is also the spectral radius of the $M\times M$
 nonnegative integer matrix $\widehat A$. If $M \ll N$, then the
 computation of $\rho(\widehat{A})$ may be feasible on a desktop
 computer whereas the computation of $\rho(A)$ may be infeasible
 on a supercomputer.

 We show that that the automorphism group of $\Theta_d(\m')$
 contains a subgroup isomorphic to the group of translations of
 $T(\m')$.  If $\Gamma_1 = \cdots = \Gamma_{d-1} = \Delta$ and
 $\Delta$ is symmetric, then the automorphism group of
 $\Theta_d(\m')$ contains a subgroup isomorphic to the group of
 rigid motions of $T(\m')$ (motions preserving the distance on
 $T(\m')$, i.e., translations, reflections and coordinate
 transpositions for equal dimensions). For example, $T(m)$ has $m$
 translations and $2m$ rigid motions if $m
 > 2$.

 We now discuss the monomer-dimer covers of $\Z^d$, see
 \cite{FoR}. A \emph{dimer} is a domino consisting of two
 neighboring atoms occupying the places $\i,\i+\e_k \in \Z^d$. A
 \emph{monomer} is a single atom occupying the place $\i \in
 \Z^d$.  A \emph{monomer-dimer cover}, respectively \emph{dimer
 cover}, of $\Z^d$ is a partition of $\Z^d$ into monomers and
 dimers, respectively dimers. We denote by $h_d$ and
 $\widetilde{h}_d$ the entropies of the monomer-dimer and dimer
 covers, respectively.  It is fairly easy to compute the values
 $h_1 = \log \frac{1+\sqrt{5}}{2}$ and $\widetilde{h}_1 = 0$. The
 big breakthrough in the sixties was a close formula for
 $\widetilde h_2$ in \cite{Fis,Kas}. The exact values of $h_d$ for
 $d \geq  2$ and $\widetilde{h}_d$ for $d \geq  3$ are unknown.

 A seminal contribution to the study of upper and lower bounds and
 estimates for $\widetilde{h}_d$ and $h_d$ was given in
 \cite{Ha1,Ha2,Ha3,HaM}. In particular, it was shown in \cite{Ha1}
 that for $p\in [0,1]$, there exists the entropy $\lambda_d(p)$ of
 the monomer-dimer covers of $\Z^d$, where $p$ is the ``density''
 of dimers, i.e., the number of dimers in the cover divided by one
 half of the volume.  The entropy $\lambda_d(p)$ is a continuous
 concave function of $p$ and $\lambda_d(1) = \widetilde{h}_d$. It
 is shown here that $h_d = \max_{p \in [0,1]} \lambda_d(p)$.  It
 was pointed out by Kingman, see \cite{Ha2}, that the van der
 Waerden conjecture for permanents of doubly-stochastic matrices
 gives a lower bound on $\widetilde{h}_d$. The improved lower
 bound for the permanents of $0$-$1$ matrices \cite{Sch} gives the
 currently best lower bound $\widetilde{h}_3 \geq 0.440075$. A
 recent breakthrough \cite{Ciu} gives the upper bound
 $\widetilde{h}_3 \leq 0.463107$, improved in \cite{Lun} to
 $\widetilde{h}_3 \leq 0.457547$.

 It is shown in \cite{Fr2} that the dimer covers can be encoded as
 $\widetilde{\Lambda}$-covers for an appropriate $d$-digraph
 $\widetilde{\Lambda} = (\widetilde{\Lambda}_1,\ldots
 ,\widetilde{\Lambda}_d)$, where all digraphs are on the set of
 vertices $\an{2d}$. We show that the monomer-dimer covers can be
 similarly encoded as $\Lambda$-configurations for an appropriate
 $d$-digraph $\Lambda = (\Lambda_1,\ldots, \Lambda_d)$, where all
 digraphs are on the set of vertices $\an{2d+1}$. Unfortunately,
 in these encodings the digraphs $\Gamma_k$,
 $\widetilde{\Gamma}_k$ are not symmetric, so (\ref{entdef}) and
 the lower bounds do not apply directly.  One of the purposes of
 this paper is to show that the entropies $h_d$ and
 $\widetilde{h}_d$ nevertheless obey upper and lower bounds
 converging to the true entropies, similar to (\ref{entdef}) and
 the lower bounds discussed above for the symmetric isotropic
 nearest neighbor digraph.
 The bounds for $h_d$ are stated in terms of the spectral radii of
 certain multidigraphs $\Theta_d(\m')$ whose automorphism group
 has a subgroup isomorphic to the the group of rigid motions of
 $T(\m')$. This fact enables us to compute the values of $h_2$ and
 $h_3$ with good precision. We also show that $\lambda_d(p)$ can
 be bounded below by using the generalized van der Waerden
 conjecture (Tverberg's conjecture), proved by the first author in
 \cite{Fr3}. For $d = 2,3$, this lower bound is better than those
 of \cite{BoW} and \cite{HaM} except for very high $p$. Our lower
 bound for $\lambda_d(p)$ yields in particular a lower bound for
 $h_d$. For $d = 2$ this lower bound is somewhat weaker than the
 one obtained from the numerical computations of
 $\rho(\Theta_d(\m'))$, but for $d = 3$ the situation is reversed.

 See \cite{HL} for a general theory of monomer-dimer covers of an
 arbitrary graph. Finally it is worth  mentioning the theoretical
 work \cite{Jer}, which shows that the general monomer-dimer
 problem in arbitrary planar graphs is computationally
 intractable.

 The contents of the paper is as follows. In Section 2 we discuss
 the the general theory of $\Z^d$ subshifts of finite type (SOFT).
 In Section 3 we prove the main inequalities of the entropy of
 $\Z^d$-SOFT with $d-1$ symmetric digraphs $\Gamma_1, \ldots,
 \Gamma_{d-1}$.  In Section 4 we recall the main features of the
 entropy of the monomer-dimer and dimer covers. In Section 5 we
 give lower bounds for the entropy of the monomer-dimer covers
 with a fixed dimer density using the lower bounds on permanents.
 In Section 6 we show that there exist analogs of the upper and
 lower bounds discussed in Section 3 that apply to the
 monomer-dimer and dimer entropy. In Section 7 we discuss using
 automorphism subgroups to reduce the computations. In Section 8
 we give numerical upper and lower bounds for $h_2$,
 $\widetilde{h}_2$, $h_3$, $\widetilde{h}_3$, and compare
 graphically our lower bounds for $\lambda_2(p)$ and $\lambda_3(p)$
 with the known lower bounds and estimates.

 \section{SOFT and NNSOFT}

 Let $\an{n}^{\Z^d}$ be the set of all colorings $\phi: \Z^d \to
 \an{n}$ of $\Z^d$ with colors from $\an{n} = \{1,\ldots,n\}$.
 Given a $d$-digraph $\Gamma = (\Gamma_1,\ldots,\Gamma_d)$ on
 $\an{n} \times \an{n}$, let $\Gamma^{\Z^d} \subseteq
 \an{n}^{\Z^d}$ be the set of all $\Gamma$-covers, namely
 colorings $\phi = (\phi_{\m})_{\m \in \Z^d}$ in $\an{n}^{\Z^d}$
 such that for each $\i \in \Z^d$ and $k \in \an{d}$, the
 restriction of $\phi$ to the line through $\i$ in the direction
 of $\e_k$, i.e., $(\phi_{\i + j\e_k})_{j \in \Z}$, is a
 bi-infinite walk on $\Gamma_k$. In ergodic theory,
 $\Gamma^{\Z^d}$ is called a \emph{nearest neighbor subshift of
 finite type (NNSOFT)}. Note that for an NNSOFT $\Gamma^{\Z^d}$
 and for $\m \in \N^d$, $W(\m)$ is the set of all configurations
 $\psi \in \an{n}^{\an{\m}}$ such that $\i, \i + \e_k \in \an{\m}$
 imply $(\psi_{\i},\psi_{\i + \e_k}) \in \Gamma_k$.

 A general SOFT can be described as follows. Let $\bM \in \N^d$
 and a nonempty subset $\cP \subseteq \an{n}^{\an{\bM}}$ be given.
 Every element $a \in \cP$ is viewed as an allowed coloring
 (configuration) of the box $\an{\bM}$ in $n$ colors. For $\i \in
 \Z^d$, we define the shifted coloring $\tau_{\i}(a)$ of $a \in
 \cP$ as the coloring of the shifted box $\an{\bM} + \i$ that
 gives to the point $\x + \i$ the same color that $a$ gives to $\x
 \in \an{\bM}$. We denote by $\tau_{\i}(\cP)$ the set
 $\{\tau_{\i}(a) : a \in \cP\}$, and regard it as the set of
 allowed colorings of $\an{\bM} + \i$. A coloring $\phi \in
 \an{n}^{\Z^d}$ is called a \emph{$\cP$-state} if for each $\i \in
 \Z^d$ the restriction of $\phi$ to $\an{\bM} + \i$ is in
 $\tau_{\i}(\cP)$. We denote by $\an{n}^{\Z^d}(\cP)$ the set of
 all $\cP$-states. In ergodic theory the set $\an{n}^{\Z^d}(\cP)$
 is called a \emph{subshift of finite type (SOFT}) \cite{Sc}.

 Each NNSOFT $\Gamma^{\Z^d}$ is a special kind of SOFT obtained by
 letting $\bM = (2,\ldots,2)$ and $\cP$ the set of all colorings
 $\psi \in \an{n}^{\an{\bM}}$ such that $\i, \i + \e_k \in \an{\bM}$
 imply $(\psi_{\i},\psi_{\i + \e_k}) \in \Gamma_k$. Conversely
 \cite{Fr1}, each SOFT $\an{n}^{\Z^d}(\cP)$ can be encoded as an
 NNSOFT $\Gamma^{\Z^d}$, where $\Gamma =
 (\Gamma_1,\ldots,\Gamma_d)$ are defined as follows. Take $N =
 \#\cP$ and use a bijection between $\cP$ and $\an{N}$. The digraph
 $\Gamma_k \subseteq \an{N} \times \an{N}$ is defined so that for
 $a,b \in \cP$ we have $(a,b) \in \Gamma_k$ if and only if there
 is a configuration $\phi \in \an{n}^{\an{\bM + \e_k}}$ such that
 the restriction of $\phi$ to $\an{\bM}$ is $a$ and the
 restriction of $\phi$ to $\an{\bM} + \e_k$ is $\tau_{\e_k}(b)$.
 Because of this equivalence, we will be dealing here with NNSOFT
 only.

 In the sequel we will be taking $\limsup$ and $\liminf$ of real
 multisequences $(a_{\m})_{\m \in \N^d}$ as $\m \to \infty$. In order
 to be clear, we define these here and prove that they are limits
 of subsequences. We also define the limit of real multisequence in
 terms of $\limsup$ and $\liminf$, which is equivalent to other
 definitions in the literature.

 \begin{defn}\label{def:multilimsup}
 Let $(a_{\m})_{\m \in \N^d}$ be a multisequence of real numbers.
 Then
 \begin{enumerate}
 \item[(a)]
 $\limsup_{\m \to \infty} a_\m$ is defined as the supremum
 (possibly $\pm \infty$) of all numbers of the form $\limsup_{q
 \to \infty} a_{\m_q}$, where $(\m_q)_{q \in \N}$ is a sequence in
 $\N^d$ satisfying $\lim_{q \to \infty} \m_q = \infty$, i.e.,
 $\lim_{q \to \infty} (\m_q)_i = \infty$ for each $i \in \an{d}$.
 We define $\liminf_{\m \to \infty} a_\m$ similarly.
 \item[(b)]
 $\lim_{\m \to \infty} a_\m = \alpha$ means
 $\limsup_{\m \to \infty} a_\m = \liminf_{\m \to \infty} a_\m =
 \alpha$.
 \end{enumerate}
 \end{defn}

 \begin{prop}\label{prop:multilimsup}
 If $\limsup_{\m \to \infty} a_\m = \alpha$, then there exists a
 sequence $(\n_q)_{q \in \N} \subseteq \N^d$ satisfying
 $\lim_{q \to \infty} \n_q = \infty$ such that the sequence
 $(a_{\n_q})_{q \in \N}$ has a limit and $\lim_{q \to \infty} a_{\n_q} =
 \alpha$. Similarly for $\liminf$.
 \end{prop}

 \proof Since the $\limsup$ of each real sequence is the limit of
 a subsequence, we may assume that we have a sequence of
 convergent subsequences $\{ a_{\m_q^i} \}$ satisfying $\lim_{q
 \to \infty} a_{\m_q^i} = \alpha_i$ and $\lim_{q \to \infty}
 \m_q^i = \infty$ for each $i \in \N$, and $\lim_{i \to \infty}
 \alpha_i = \alpha$. Note that $\alpha_i \leq \alpha$ for all $i$
 by definition of the supremum. If $\alpha_i = \alpha$ for some
 $i$, we are done. In particular, if $\alpha = -\infty$, then
 $\alpha_i = -\infty = \alpha$ for each $i$ and we are done.
 Therefore we may assume that $\alpha \in \R \cup \{\infty\}$ and
 that $\alpha_i < \alpha_{i+1}$ for all $i$.

 Assume first that $\alpha \in \R$. Then for each $i \in \N$ there
 exists a $q(i) \in \N$ such that $\m_{q(i+1)}^{i+1} > 2
 \m_{q(i)}^i$ and $|a_{\m_{q(i)}^i} - \alpha_i| < \frac{1}{2^i}$.
 Then we can take $\n_i = \m_{q(i)}^i$ and the result follows.
 Similarly, if $\alpha = \infty$, then for each $i \in \N$ there
 exists a $q(i) \in \N$ such that $\m_{q(i+1)}^{i+1} > 2
 \m_{q(i)}^i$ and $a_{\m_{q(i)}^i} > \alpha_i - 1$, and again we can take
 $\n_i = \m_{q(i)}^i$.
 \qed

 Let $W_\perio(\m) \subseteq \Gamma^{\Z^d}$ be the
 set of periodic $\Gamma$-covers with period $\m$.  Then
 \begin{equation}\label{perent}
 h_{\perio}(\Gamma) := \limsup_{\m \to \infty} \frac{\log \#W_{\perio}
 (\m)}{|\m|_\pr}
 \end{equation}
 is called the \emph{periodic entropy} of $\Gamma^{\Z^d}$.
 Clearly $ h_{\perio}(\Gamma) \leq  h(\Gamma)$.

 \section{Main Inequalities for Symmetric NNSOFT}\label{sec:SymmetricNNSOFT}

 For $d \geq 2$, consider $\m =
 (m_1,\ldots,m_d) \in \N^d$ and $\m^{-} := (m_2,\ldots,m_d)$. Let
 $W_{\perio,\{1\}}(\m)$ be the set of $\Gamma$-configurations
 in the box $\an{\m}$ that correspond to $\Gamma$-covers of
 $T(m_1) \times \an{\m^{-}}$, i.e., that can be extended
 periodically in the direction of $\e_1$ with period $m_1$ into
 $\Gamma$-covers of $\Z \times \an{\m^{-}}$. We can view these
 configurations as $\widehat{\Gamma}$-configurations in the box
 $\an{\m^{-}}$, where $\widehat{\Gamma} =
 (\widehat{\Gamma}_2,\ldots,\widehat{\Gamma}_d)$, for each $j$ the
 vertex set of $\widehat{\Gamma}_j$ is the set
 $\Gamma_{1,\perio}^{m_1}$ of closed walks $a =
 (a_1,\ldots,a_{m_1},a_1)$ of length $m_1$ on $\Gamma_1$, and
 where $(a,b) \in \widehat{\Gamma}_j$ if and only if $(a_i,b_i)
 \in \Gamma_j$ for $i =1,\ldots,m_1$. For this reason, the
 following limit exists and is equal to the entropy
 $h(\widehat{\Gamma}_2,\ldots,\widehat{\Gamma}_d)$ of the NNSOFT
 $\widehat{\Gamma}^{{\Z}^{d-1}}$:
 \begin{equation}\label{defbarh}
  \overline{h}(m_1,\Gamma) := \lim_{\m^- \to \infty} \frac{\log
  \#W_{\perio,\{1\}}(\m)}{|\m^{-}|_{\pr}}, \quad m_1 \in \N \quad
 \end{equation}
 We define $W^{-}(\m^{-})$ as the set of
 $(\Gamma_2,\ldots,\Gamma_d)$-covers of the box $\an{\m^{-}}$. In
 the degenerate case $m_1=0$, we define
 $W_{\perio,\{1\}}(0,\m^{-})$ to be simply $W^{-}(\m^{-})$ and
 $(\widehat{\Gamma}_2,\ldots,\widehat{\Gamma}_d)$ to be simply
 $(\Gamma_2,\ldots,\Gamma_d)$. Then (\ref{defbarh}) is also valid
 for $m_1=0$, where we understand $\overline{h}(0,\Gamma)$ to be
 $h(\Gamma_2,\ldots ,\Gamma_d)$.

 \begin{theo}\label{perulb}  Consider the NNSOFT
  $\;\Gamma^{\Z^d}$ for $d \geq 2$.
  Let $h(\Gamma)$ and $\overline{h}(r,\Gamma)$ be defined by (\ref{hGamma}) and
  (\ref{defbarh}), respectively.
  Assume that $\Gamma_1$ is symmetric.  Then for all $p,r \in \N$ and $q \in \Z_+$,
  \begin{equation}\label{perulb1}
   \frac{\overline{h}(2r,\Gamma)}{2r}  \geq  h(\Gamma)  \geq
   \frac{\overline{h}(p + 2q,\Gamma)- \overline{h}(2q,\Gamma)}{p}.
  \end{equation}
 \end{theo}
 \proof  Fix $\m^{-} = (m_2,\ldots,m_d) \in \N^{d-1}$ and let
 $\Omega_1(\m^{-})$ be the following transfer digraph on the
 vertex set $W^{-}(\m^{-})$, analogous to the transfer digraph
 $\Omega_d(\m')$ described in Section 1. Vertices $u,v$ satisfy $(u,v)
 \in \Omega_1(\m^{-})$ if and only if $[u,v] \in W(2,\m^{-})$,
 where $[u,v]$  is the configuration consisting of $u,v$ occupying
 the levels $x_1=1,2$ of $\an{(2,\m^{-})}$, respectively. Let
 $N=\# W^{-}(\m^{-})$ and let $C(\m^{-})$ be the $N \times N$
 $0$-$1$ incidence matrix of $\Omega_1(\m^{-})$, with spectral
 radius $\rho(C(\m^{-})$. As a nonnegative matrix, $C(\m^{-})$
 satisfies (see e.g., \cite{Fr2})
 \[\log \rho(C(\m^{-})) = \lim_{k \to \infty} \frac{\log \1\trans C(\m^{-})^k \1}{k},\]
 where $\1 = (1,\ldots,1)\trans$. Since $\1\trans C(\m^{-})^k \1$
 is the number of walks of length $k$ on $\Omega_1(\m^{-})$, which
 correspond to $\Gamma$-covers of $\an{(k,\m^{-})}$, we obtain
 \begin{equation}\label{logrhoeqlim}
 \frac{\log\rho(C(\m^{-}))} {|\m^{-}|_{\pr}} =
 \lim_{k\to\infty} \frac{\log\#W(k,\m^{-})}{k\,|\m^{-}|_{\pr}}.
 \end{equation}
 Now send $m_2,\ldots,m_d$ to $\infty$, and observe that by
 (\ref{hGamma}) and (\ref{ub1}), the right-hand side of
 (\ref{logrhoeqlim}) converges to $h(\Gamma)$ and is an upper
 bound on it for each $\m^{-}$. Thus we obtain \cite{Fr1}
 \begin{align}
 \frac{\log\rho(C(\m^{-}))}{|\m^{-}|_{\pr}} &\geq
 h(\Gamma), \qquad \m^{-} \in \N^{d-1}\label{enteq}\\
 \lim_{\m^- \to \infty}\frac{\log\rho(C(\m^{-}))}
 {|\m^{-}|_{\pr}} &= h(\Gamma).\label{enteq1}
 \end{align}
 Next, we observe that
 \begin{equation}\label{percount}
 \tr C(\m^{-})^{q} = \#W_{\perio,\{1\}}(q,\m^{-}), \qquad q \in \Z_+,
 \end{equation}
 where $C(\m^{-})^0$ is the $N \times N$ identity matrix. Recall
 that the trace of $C(\m^{-})^q$ is given by
 \[ \tr C(\m^{-})^{q}=\sum_{i=1}^N \lambda_i^{q}, \qquad q \in \Z_+,\]
 where $\lambda_1,\ldots,\lambda_N$ be the eigenvalues of
 $C(\m^{-})$. Since $C(\m^{-})$ is a nonnegative matrix, the
 Perron-Frobenius theorem yields that its spectral radius
 $\rho(C(\m^{-})):=\max_{i \in \an{N}} |\lambda_i|$ is one of the
 $\lambda_i$. Since by assumption $\Gamma_1$ is symmetric,
 $\Omega_1(\m^{-})$ and hence $C(\m^{-})$ are symmetric. Therefore
 $\lambda_1, \ldots,\lambda_N$ are real, and hence $\tr
 C(\m^{-})^{2r} \geq  \rho(C(\m^{-}))^{2r}$ for each $r \in \N$.
 Taking logarithms and using (\ref{percount}), we obtain
 \begin{equation}\label{basineq}
 \frac{\log \#W_{\perio,\{1\}}(2r,\m^{-})}{2r |\m^{-}|_{\pr}} \geq
 \frac{\log\rho(C(\m^{-}))}{|\m^{-}|_{\pr}}, \qquad r \in \N.
 \end{equation}
 Sending $m_2,\ldots,m_d$ to $\infty$ in (\ref{basineq}) and using
 (\ref{defbarh}) and (\ref{enteq1}), we
 deduce the upper bound on $h(\Gamma)$ in (\ref{perulb1}).

 To prove the lower bound in (\ref{perulb1}), we note that
 \begin{multline*}
 \tr C(\m^{-})^{p+2q} = \sum_i \lambda_i^{p+2q} \leq \sum_i |\lambda_i|^{p+2q}
 = \sum_i |\lambda_i|^p \lambda_i^{2q}\\
 \leq \sum_i \rho(C(\m^{-}))^p \lambda_i^{2q} = \rho(C(\m^{-}))^p \tr C(\m^{-})^{2q}
 \end{multline*}
 and thus by (\ref{percount})
 \begin{gather}
 \rho(C(\m^{-}))^p  \geq  \frac{\tr C(\m^{-})^{p+2q}}{\tr
 C(\m^{-})^{2q}}=\frac{\#W_{\perio,\{1\}}(p+2q,\m^{-})}
 {\#W_{\perio,\{1\}}(2q,\m^{-})}\label{rhogeqtraceratio}
 \\
 \frac{\log\rho(C(\m^{-}))}{|\m^{-}|_{\pr}}  \geq
 \frac{\log\#W_{\perio,\{1\}}(p+2q,\m^{-}) -
 \log \#W_{\perio,\{1\}}(2q,\m^{-})}{p|\m^{-}|_{\pr}}.\notag
 \end{gather}
 Sending $\m^-$ to $\infty$ and using (\ref{enteq1}) and
 (\ref{defbarh}) (recall that the latter holds for $m_1 \in \Z_+$), we
 deduce the lower bound in (\ref{perulb1}).  \qed

 When $d=2$, $\overline{h}(m_1,\Gamma)$ is the entropy of the NNSOFT
 $\widehat{\Gamma}_2^\Z$ (recall that $\widehat{\Gamma}_2$ is
 simply $\Gamma_2$ when $m_1 = 0$). Since this is a
 $1$-dimensional NNSOFT, that entropy is equal to
 $\log\rho(\widehat{\Gamma}_2)$. We denote
 $\rho(\widehat{\Gamma}_2)$ by $\theta_2(m_1)$, and obtain the
 following corollary to Theorem~\ref{perulb}.

 \begin{corol}\label{ulbd=2}
 Let $d=2$ and assume that $\Gamma_1$ is symmetric.
 Then for all $p,r \in \N$ and $q \in \Z_+$,
 \begin{equation}\label{ubspecr2}
 \frac{\log\theta_2(2r)}{2r} \geq  h(\Gamma) \geq
 \frac{\log\theta_2(p+2q)-\log\theta_2(2q)}{p},
 \end{equation}
 where $\theta_2$ is defined above.
 \end{corol}

 In (\ref{ubspecr2}) take $q=0$ and $p=2r$, and send $r$ to
 $\infty$.  Clearly the upper and lower bounds then converge to
 $h(\Gamma)$. Hence $h(\Gamma)$ is computable \cite{Fr1}. For
 completeness of the exposition we reproduce a short proof of
 (\ref{entdef}) for any $d \geq  2$ given in \cite{Fr2}. We use
 the following straightforward lemma.

 \begin{lemma}\label{TransferGammaHat}
 Let $\Gamma = (\Gamma_1,\ldots,\Gamma_d)$ and $\m \in \N^d$, put
 $\Gamma' = (\Gamma_1,\ldots,\Gamma_{d-1})$ and $\m' =
 (m_1,\ldots,m_{d-1})$, and let $\Theta_d(\m')$ be the transfer
 digraph between $\Gamma'$-covers of $T(\m')$ with respect to
 $\Gamma_d$. Let $\widehat{\Gamma}_2,\dots,\widehat{\Gamma}_d$ be
 defined as in the beginning of this section, put
 $\widehat{\Gamma}\,' =
 (\widehat{\Gamma}_2,\ldots,\widehat{\Gamma}_{d-1})$ and
 $\widetilde{\m} = (m_2,\ldots,m_{d-1})$, and let
 $\widehat{\Theta}_d(\widetilde{\m})$ be the transfer digraph
 between $\widehat{\Gamma}\,'$-covers of $T(\widetilde{\m})$ with
 respect to $\widehat{\Gamma}_d$. Then $\Theta_d(\m')$ and
 $\widehat{\Theta}_d(\widetilde{\m})$ are isomorphic, and in
 particular $\rho(\Theta_d(\m')) =
 \rho(\widehat{\Theta}_d(\widetilde{\m}))$.
 \end{lemma}
 \proof We use the following bijection between the vertices $u$ of
 $\Theta_d(\m')$ and the vertices $\widehat{u}$ of
 $\widehat{\Theta}_d(\widetilde{\m})$. Given $u = (\phi_{\i})_{\i \in
 \an{\m'}}$, we have
 \begin{equation}\label{AdjInGammak}
  (\phi_{\i},\phi_{\i + \e_k}) \in \Gamma_k, \quad k = 1,\ldots,d-1,
 \end{equation}
 where the addition $\i + \e_k$ is understood modulo $m_k$, i.e.,
 $m_k + 1$ is $1$. Then the corresponding $\widehat{u}$ is defined
 to be $\widehat{u} = (\widehat{\phi}_\j)_{\j \in
 \an{\widetilde{\m}}}$, where $\widehat{\phi}_{\j} =
 (\phi_{(q,\j)})_{q=1}^{m_1}$. We note that $\widehat{\phi}_{\j}$
 is indeed a $\Gamma_1$-cover of $T(m_1)$ and thus a vertex of
 $\widehat{\Gamma}_2,\ldots,\widehat{\Gamma}_{d-1}$ by
 (\ref{AdjInGammak}) with $\i = (q,\j)$ and $k = 1$. In order to
 show that $\widehat{u}$ is a $\widehat{\Gamma}\,'$-cover of
 $T(\widetilde{\m})$ and thus a vertex of
 $\widehat{\Theta}_d(\widetilde{\m})$, we need to show that
 $(\widehat{\phi}_{\j},\widehat{\phi}_{\j + \e'_k}) \in
 \widehat{\Gamma}_k$ for $k = 2,\ldots,d-1$. This means showing
 that $(\phi_{(q,\j)},\phi_{(q,\j+\e'_k)}) \in \Gamma_k$ for $k =
 2,\ldots,d-1$ and $q=1,\ldots,m_1$, which follows in turn from
 (\ref{AdjInGammak}) with $\i = (q,\j)$. It is easy to see that
 the correspondence $u \mapsto \widehat{u}$ can be inverted. It
 remains to show that $(u,v) \in \Theta_d(\m') \Leftrightarrow
 (\widehat{u},\widehat{v}) \in
 \widehat{\Theta}_d(\widetilde{\m})$. We prove only the
 $\Rightarrow$ part. Let $u = (\phi_{\i})_{\i \in \an{\m'}}$ and
 $v = (\psi_{\i})_{\i \in \an{\m'}}$ be $\Gamma'$-covers of
 $T(\m')$. The assumption $(u,v) \in \Theta_d(\m')$ means that
 $(\phi_\i,\psi_\i) \in \Gamma_d$ for all $\i \in \an{\m'}$.
 Applying this with $\i = (q,\j)$, $q = 1,\ldots,m_1$ and $\j \in
 \an{\widetilde{\m}}$ shows that
 $(\widehat{\phi}_\j,\widehat{\psi}_\j) \in \widehat{\Gamma}_d$
 for all $\j \in \an{\widetilde{\m}}$, which means in turn that
 $(\widehat{u},\widehat{v}) \in
 \widehat{\Theta}_d(\widetilde{\m})$.
 \qed
 \begin{theo}\label{ubspecrt}  Let $d \geq  2$ and consider the
 NNSOFT $\Gamma^{\Z^d}$, where $\Gamma =
 (\Gamma_1,\ldots,\Gamma_{d})$. For $\m'=(m_1,\ldots ,m_{d-1}) \in
 \N^{d-1}$ and $\Gamma'=(\Gamma_1,\ldots ,\Gamma_{d-1})$, let
 $\Theta_d(\m')$ be the transfer digraph between $\Gamma'$-covers
 of $T(\m')$ with respect to $\Gamma_d$. Assume that
 $\Gamma_1,\ldots ,\Gamma_{d-1}$ are symmetric and $m_1,\ldots
 ,m_{d-1}$ are even. Then
 \[h(\Gamma) \leq \frac{\log\rho(\Theta_d(\m'))}{|\m'|_{\pr}}.\]
 \end{theo}
 \proof The proof is by induction on $d$.  For $d = 2$ the result
 is equivalent to the upper bound in (\ref{ubspecr2}). For the
 induction step, observe that the upper bound of (\ref{perulb1})
 with $r = m_1/2$ yields $h(\Gamma) \leq \overline{h}(m_1,\Gamma)/m_1$.
 Recall that $\overline{h}(m_1,\Gamma)$ is the entropy of the NNSOFT
 $\widehat{\Gamma}^{\Z^{d-1}}$, where $\widehat{\Gamma} =
 (\widehat{\Gamma}_2,\ldots,\widehat{\Gamma}_d)$ is as in
 Lemma~\ref{TransferGammaHat}. Since
 $\Gamma_2,\ldots,\Gamma_{d-1}$ are symmetric, so are
 $\widehat{\Gamma}_2,\ldots,\widehat{\Gamma}_{d-1}$, and therefore
 the induction hypothesis applied to $\widehat{\Gamma}^{\Z^{d-1}}$
 gives $\overline{h}(m_1,\Gamma) \leq
 \log\rho(\widehat{\Theta}_d(\widetilde{\m}))/|\widetilde{\m}|_{\pr}$,
 where $\widetilde{\m}$ and $\widehat{\Theta}_d$ are as in the
 lemma. Finally, an application of the lemma completes the proof.
  \qed
 \begin{corol}\label{lubd=3}
 Let $\Gamma=(\Gamma_1,\Gamma_2,\Gamma_3)$ and assume that
 $\Gamma_1$ and $\Gamma_2$ are symmetric. For $(m_1,m_2) \in
 \N^2$, let $\Theta_3(m_1,m_2)$ be the transfer digraph between
 $(\Gamma_1,\Gamma_2)$-covers of $T(m_1,m_2)$ with respect to
 $\Gamma_3$, and let $\theta_3(m_1,m_2)$ be its spectral radius.
 Let $\theta_3(0,m_2)$ be the spectral radius of the transfer
 digraph between $\Gamma_2$-covers of $T(m_2)$ with respect to
 $\Gamma_3$. Let $\theta_3(m_1,0)$ be the spectral radius of the
 transfer digraph between $\Gamma_1$-covers of $T(m_1)$ with
 respect to $\Gamma_3$. Then for all $r,t,p,u,v \in \N$ and $q,s
 \in \Z_+$ we have
 \begin{multline}\label{lubd=3i}
 \frac{\log\theta_3(2r,2t)}{4rt} \geq  h(\Gamma)\\
 \geq \frac{\log\theta_3(p+2q,u+2s) - \log\theta_3(p+2q,2s)}{up} -
 \frac{\log\theta_3(2q,2v)}{2vp}.
 \end{multline}
 \end{corol}
 \proof   The upper bound in (\ref{lubd=3i}) follows directly from
 Theorem~\ref{ubspecrt} for $d=3$. To show the lower bound we use
 the lower bound in (\ref{perulb1}), which is valid since
 $\Gamma_1$ is symmetric, and gives
 \begin{equation}\label{lb1}
 h(\Gamma_1,\Gamma_2,\Gamma_3) \geq
 \frac{\overline{h}(p+2q,(\Gamma_1,\Gamma_2,\Gamma_3)) -
 \overline{h}(2q,(\Gamma_1,\Gamma_2,\Gamma_3))}{p}.
 \end{equation}
 For each $a \in \N$ we have
 $\overline{h}(a,(\Gamma_1,\Gamma_2,\Gamma_3)) =
 h(\widehat{\Gamma}_2,\widehat{\Gamma}_3)$, where
 $\widehat{\Gamma}_2,\widehat{\Gamma}_3$ are digraphs on the
 vertex set $\Gamma_{1,\perio}^a$ as in the beginning of this
 section. Since $\Gamma_2$ is symmetric, so is
 $\widehat{\Gamma}_2$, and so we can apply the lower bound of
 Corollary \ref{ulbd=2} to
 $(\widehat{\Gamma}_2,\widehat{\Gamma}_3)$ to obtain
 \begin{equation}\label{lb2}
 h(\widehat{\Gamma}_2,\widehat{\Gamma}_3) \geq
 \frac{\log\theta_3(a,u+2s) - \log\theta_3(a,2s)}{u},
 \end{equation}
 where $\theta_3(a,b) = \rho(\widehat{\Theta}_3(b)) =
 \rho(\Theta_3(a,b))$ by Lemma \ref{TransferGammaHat}. Inequality
 (\ref{lb2}) is also valid for $s=0$, since we defined
 $\theta_3(a,0)$ to be the spectral radius of
 $\widehat{\Gamma}_3$, exactly as in Corollary \ref{ulbd=2} for the
 degenerate case. Using (\ref{lb2}) for $a=p + 2q$ gives
 \begin{equation}\label{lb3}
 \overline{h}(p+2q,(\Gamma_1,\Gamma_2,\Gamma_3)) \geq
 \frac{\log\theta_3(p+2q,u+2s) - \log\theta_3(p+2q,2s)}{u}.
 \end{equation}
 Apply the upper bound of Corollary \ref{ulbd=2} to
 $(\widehat{\Gamma}_2,\widehat{\Gamma}_3)$ to obtain
 \begin{equation}\label{lb4}
 h(\widehat{\Gamma}_2,\widehat{\Gamma}_3) \leq
 \frac{\log\theta_3(a,2v)}{2v}.
 \end{equation}
 Inequality (\ref{lb4}) is also valid for $a=0$, since in that case
 $(\widehat{\Gamma}_2,\widehat{\Gamma}_3) = (\Gamma_2,\Gamma_3)$,
 by Theorem \ref{ubspecrt} applied to $(\Gamma_2,\Gamma_3)$, and
 by the definition of $\theta_3(0,2v)$. Using (\ref{lb4}) for
 $a=2q$ gives
 \begin{equation}\label{lb5}
 \overline{h}(2q,(\Gamma_1,\Gamma_2,\Gamma_3)) \leq
 \frac{\log\theta_3(2q,2v)}{2v}.
 \end{equation}
 Finally, substitution of (\ref{lb3}) and (\ref{lb5}) in
 (\ref{lb1}) yields the lower bound of (\ref{lubd=3i}). \qed

 \section{Dimer and Monomer-Dimer Covers of $\Z^d$}

 As in \cite{Fr2}, the set of monomer-dimer covers, respectively
 dimer covers, of $\Z^d$ is an NNSOFT $\Gamma^{\Z^d}$,
 respectively $\widetilde{\Gamma}^{\Z^d}$, where $\Gamma$ and
 $\widetilde{\Gamma}$ are defined as follows. We encode a
 monomer-dimer cover of $\Z^d$ as a coloring of $\Z^d$ with the
 $2d + 1$ colors $1,\ldots,2d + 1$: a dimer in the direction of
 $\e_k$ occupying the adjacent points $\i,\i + \e_k$ is encoded by
 the color $k$ at $\i$ and the color $k + d$ at $\i + \e_k$; a
 monomer at $\i$ is encoded by the color $2d + 1$ at $\i$.  This
 imposes restrictions on the coloring, which are expressed by the
 $d$-digraph $\Gamma = (\Gamma_1,\ldots,\Gamma_k)$ on the set of
 vertices $\an{2d + 1}$, where
 \begin{itemize}
  \item $(k,q) \in \Gamma_k \Leftrightarrow q = k + d$;
  \item for $j \neq k$, $(j,q) \in \Gamma_k \Leftrightarrow q \neq k +
  d$.
 \end{itemize}
 It is easy to check that this gives a bijection between the
 monomer-dimer covers of $\Z^d$ and $\Gamma^{\Z^d}$. Similarly, if
 $\widetilde{\Gamma} =
 (\widetilde{\Gamma}_1,\ldots,\widetilde{\Gamma}_d)$ is obtained
 from $\Gamma$ by removing the vertex $2d +1$, then there is a
 bijection between the dimer covers of $\Z^d$ and
 $\widetilde{\Gamma}^{\Z^d}$.

 The disadvantage of these encodings is that $\Gamma_k$ and
 $\widetilde{\Gamma}_k$ are not symmetric, so we cannot apply the
 results of Section \ref{sec:SymmetricNNSOFT} directly. However,
 as pointed out in \cite{Ciu} for the dimer problem, there is a
 hidden symmetry, which enables us to obtain results analogous to
 those of Section \ref{sec:SymmetricNNSOFT}.

\begin{figure}
   \begin{center}
   \psset{unit=0.5cm}
   \begin{pspicture}(0,-1)(24,8)
   \psset{linewidth=0.4pt}
   \psline[linewidth=0.5mm]{->}(1,1)(5,1)
   \uput[0](5,1){$1$}
   \psline[linewidth=0.5mm]{->}(1,1)(1,5)
   \uput[90](1,5){$2$}
   \psline(1,1)(4,1)
   \psline(1,2)(4,2)
   \psline(1,3)(4,3)
   \psline(1,4)(4,4)
   \psline(1,1)(1,4)
   \psline(2,1)(2,4)
   \psline(3,1)(3,4)
   \psline(4,1)(4,4)
   \rput(1.5,1.5){$5$}
   \rput(2.5,1.5){$4$}
   \rput(3.5,1.5){$5$}
   \rput(1.5,2.5){$5$}
   \rput(2.5,2.5){$2$}
   \rput(3.5,2.5){$1$}
   \rput(1.5,3.5){$3$}
   \rput(2.5,3.5){$4$}
   \rput(3.5,3.5){$5$}
   \rput(2.5,-1){(a)}
   \psline[linewidth=0.5mm]{->}(9,1)(14,1)
   \uput[0](14,1){$1$}
   \psline[linewidth=0.5mm]{->}(9,1)(9,6)
   \uput[90](9,6){$2$}
   \psline(8,0)(13,0)
   \psline(8,1)(13,1)
   \psline(8,2)(13,2)
   \psline(8,3)(13,3)
   \psline(8,4)(13,4)
   \psline(8,5)(13,5)
   \psline(8,0)(8,5)
   \psline(9,0)(9,5)
   \psline(10,0)(10,5)
   \psline(11,0)(11,5)
   \psline(12,0)(12,5)
   \psline(13,0)(13,5)
   \psset{linewidth=0.15,dimen=outer,linecolor=lightgray,fillstyle=solid,fillcolor=lightgray}
   \psframe(9.1,1.1)(9.9,1.9)
   \psframe(9.1,2.1)(9.9,2.9)
   \psframe(8.1,3.1)(9.9,3.9)
   \psframe(10.1,0.1)(10.9,1.9)
   \psframe(10.1,2.1)(10.9,3.9)
   \psframe(11.1,1.1)(11.9,1.9)
   \psframe(11.1,2.1)(12.9,2.9)
   \psframe(11.1,3.1)(11.9,3.9)
   \rput(10.5,-1){(b)}
   \psset{linecolor=black}
   \psline[linewidth=0.5mm]{->}(17,1)(23,1)
   \uput[0](23,1){$1$}
   \psline[linewidth=0.5mm]{->}(17,1)(17,7)
   \uput[90](17,7){$2$}
   \psset{linewidth=0.4pt}
   \psline(17,1)(22,1)
   \psline(17,2)(22,2)
   \psline(17,3)(22,3)
   \psline(17,4)(22,4)
   \psline(17,5)(22,5)
   \psline(17,6)(22,6)
   \psline(17,1)(17,6)
   \psline(18,1)(18,6)
   \psline(19,1)(19,6)
   \psline(20,1)(20,6)
   \psline(21,1)(21,6)
   \psline(22,1)(22,6)
   \psset{linewidth=0.15,dimen=outer,linecolor=lightgray}
   \psframe(18.1,2.1)(18.9,2.9)
   \psframe(18.1,3.1)(18.9,3.9)
   \psframe(17.1,4.1)(18.9,4.9)
   \psframe(19.1,1.1)(19.9,2.9)
   \psframe(19.1,3.1)(19.9,4.9)
   \psframe(20.1,2.1)(20.9,2.9)
   \psframe(20.1,3.1)(21.9,3.9)
   \psframe(20.1,4.1)(20.9,4.9)
   \rput(19.5,-1){(c)}
   \end{pspicture}
  \end{center}
  \caption{(a) $\Gamma$-coloring of $\an{\m} = \an{(3,3)}$; (b)
  Corresponding monomer-dimer cover of $T$; (c) Corresponding
  monomer-dimer cover of $S$}
  \label{protrusions}
 \end{figure}
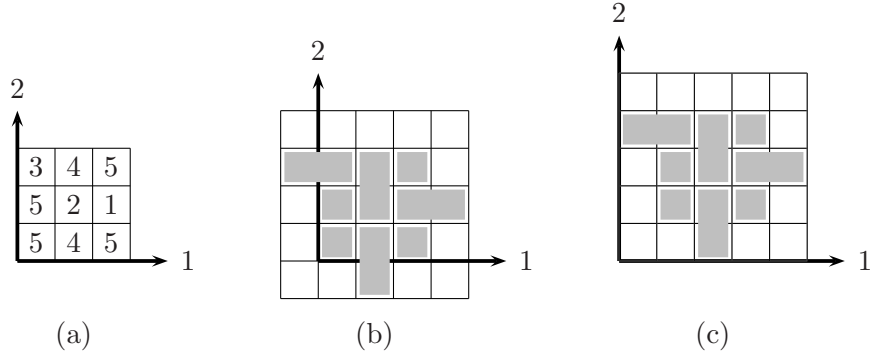

 Recall that $W(\m)$ denotes the set of $\Gamma$-colorings of
 $\an{\m} \subseteq \N^d$. Consider a $\Gamma$-coloring $\phi \in
 W(\m)$ with the $\Gamma$ defined above. Certain points $\i$ on
 the boundary of $\an{\m}$ can receive colors indicating that $\i$
 is one half of a dimer whose other half is outside $\an{\m}$.
 Therefore $\phi$ corresponds to a monomer-dimer cover of a ``box
 with protrusions'' $T$ satisfying $\an{\m} \subseteq T \subseteq
 \an{\m + 2\1} - \1$, where $\1 := (1,\ldots,1) \in \N^d$, such
 that each monomer in the cover is contained in $\an{\m}$ and each
 dimer in the cover has a nonempty intersection with $\an{\m}$. We
 translate $T$ by $\1$ to move it into $\N^d$, and thus $\phi$
 corresponds to a monomer-dimer cover of a set $S$ satisfying
 $\an{\m} + \1 \subseteq S \subseteq \an{\m + 2\1}$ such that each
 monomer in the cover is contained in $\an{\m} + \1$ and each
 dimer in the cover has a nonempty intersection with $\an{\m} +
 \1$. Conversely, each monomer-dimer cover of such a set $S$
 satisfying these conditions corresponds to a $\Gamma$-coloring of
 $\an{\m}$.  This is illustrated in Figure \ref{protrusions}.

 Similarly, $\widetilde{W}(\m)$ denotes the set of
 $\widetilde{\Gamma}$-colorings of $\an{\m}$, and there is a
 bijection between $\widetilde{W}(\m)$ and the set of dimer covers
 of a set $S$ satisfying $\an{\m} + \1 \subseteq S \subseteq
 \an{\m + 2\1}$ such that each dimer in the cover has a nonempty
 intersection with $\an{\m} + \1$.

 Let $W_{\perio}(\m)$, respectively
 $\widetilde{W}_{\perio}(\m)$, denote the set of
 $\Gamma$-colorings, respectively $\widetilde\Gamma$-colorings, of
 $\an{\m}$ that can be extended periodically to
 $\Gamma$-colorings, respectively $\widetilde{\Gamma}$-colorings,
 of $\Z^d$ with period $\m$. It corresponds to the set of
 monomer-dimer covers, respectively dimer covers, of $T(\m)$ and
 satisfies $W_{\perio}(\m) \subseteq W(\m)$,
 $\widetilde{W}_{\perio}(\m) \subseteq \widetilde{W}(\m)$.

 Finally, let $W_0(\m)$, respectively $\widetilde W_0(\m)$, be the
 set of $\Gamma$-colorings of $\an{\m}$ for which $S$ defined
 above is equal to $\an{\m} + \1$, i.e., each dimer in the
 corresponding cover of $S$ is contained in $\an{\m}$. To
 emphasize the fact that the dimers do not protrude out of
 $\an{\m}$, we refer to these covers as \emph{tilings}. We have
 $W_0(\m)\subseteq W_{\perio}(\m)$, $\widetilde W_0(\m)\subseteq
 \widetilde W_{\perio}(\m)$. We can see that $\# W(\m) \leq \#
 W_0(\m + 2\1)$, because we can extend the monomer-dimer cover of
 $S$ into a member of $W_0(\an{\m + 2\1})$ by tiling $\an{\m +
 2\1} \setminus S$ with monomers.

 From the discussion above we have
 \begin{gather}
  \# W_0(\m) \leq  \# W_{\perio}(\m) \leq  \# W(\m)  \leq  \# W_0(\m + 2\1)\label{wina}\\
  \# \widetilde{W}_0(\m) \leq  \# \widetilde{W}_{\perio}(\m) \leq  \# \widetilde{W}(\m)\label{winb}\\
  \# \widetilde{W}_0(\m) \leq  \# W_0(\m)\label{winc}\\
  \# \widetilde{W}_{\perio}(\m) \leq  \# W_{\perio}(\m)\label{wind}\\
  \#\widetilde{W}(\m) \leq  \# W(\m).\label{wine}
 \end{gather}

 Recall that the $d$-dimensional monomer-dimer
 entropy $h_d$ is defined by
 \[h_d := \lim_{\m \to \infty} \frac{\log \# W(\m)}{|\m|_{\pr}}.\]
 From (\ref{wina}) we obtain
 \begin{multline*}
 \liminf_{\m \to \infty} \frac{\log \# W_0(\m)}{|\m|_{\pr}}
 = \liminf_{\m \to \infty} \frac{\log \# W_0(\m + 2\1)}{|\m + 2\1|_{\pr}}
 \\
 = \liminf_{\m \to \infty} \frac{\log \# W_0(\m + 2\1)}{|\m|_{\pr}}
 \geq h_d \geq
 \limsup_{\m \to \infty} \frac{\log \#
 W_0(\m)}{|\m|_{\pr}}.
 \end{multline*}
 This and one more application of (\ref{wina}) give
 \begin{multline}\label{dmen}
 h_d := \lim_{\m \to \infty} \frac{\log \# W(\m)}{|\m|_{\pr}} =
 \lim_{\m \to \infty} \frac{\log \# W_{\perio}(\m)}{|\m|_{\pr}}
 \\
 = \lim_{\m \to \infty} \frac{\log \# W_0(\m)}{|\m|_{\pr}}.
 \end{multline}
 Similarly, the $d$-dimensional dimer entropy $\widetilde{h}_d$ is
 defined by
 \[\widetilde{h}_d := \lim_{\m \to \infty} \frac{\log \#
 \widetilde{W}(\m)}{|\m|_{\pr}}.\]
 It is known to satisfy
 \begin{multline}\label{den}
 \widetilde{h}_d := \lim_{\m \to \infty} \frac{\log \# \widetilde W(\m)}{|\m|_{\pr}}=
 \lim_{\m \to \infty, \frac{|\m|_{\pr}}{2} \in \N}
 \frac{\log \# \widetilde W_{\perio}(\m)}{|\m|_{\pr}}\\
 = \lim_{\m \to \infty, \frac{|\m|_{\pr}}{2}\in \N}
 \frac{\log \# \widetilde W_0(\m)}{|\m|_{\pr}}.
 \end{multline}
 The proof of (\ref{den}) is more involved, and follows from the
 results proved in \cite{Ha1}, as we show now. For $\m \in \N^d$
 and $s \in \left[0,\frac{|\m|_{\pr}}{2}\right] \cap \Z$, let
 $W_0(\m,s)$ be the subset of $W_0(\m)$ consisting of the
 monomer-dimer tilings of $\an{\m}$ that have exactly $s$ dimers.
 As pointed out in \cite{Ha1}, $W_0(\m,s) \ne \emptyset$ by
 induction on $d$.  It is shown in \cite{Ha1} that there exists a
 function $\lambda_d(\cdot) : [0,1] \to \R_+$ such that for all
 sequences $(\m_q)_{q \in \N}$ and $(s_q)_{q \in \N}$ satisfying
 \begin{equation}\label{HammCond}
 s_q \in \left[0,\frac{|\m_q|_{\pr}}{2}\right] \cap \Z, \quad
 \lim_{q \to \infty} \m_q =\infty, \quad
 \lim_{q \to \infty} \frac{2s_q}{|\m_q|_{\pr}} = p \in [0,1],
 \end{equation}
 the following equality holds
 \begin{equation}\label{lambpdef}
 \lim_{q \to \infty} \frac{\log \#W_0(\m_q,s_q)}{|\m_q|_{\pr}} = \lambda_d(p).
 \end{equation}
 Furthermore, the function $\lambda_d(p)$ is a continuous concave function of $p$ on $[0,1]$.
 We call $\lambda_d(p)$ the \emph{monomer-dimer entropy with dimer
 density $p$}.
 \begin{theo}\label{maxlambp}
 Let $\widetilde{W}(\m)$, $\widetilde{W}_\perio(\m)$,
 $\widetilde{W}_0(\m)$ be defined as above. Then (\ref{den}) and
 the following equalities hold
 \begin{gather}
 \lambda_d(0) = 0 \label{lambp1}\\
 \lambda_d(1) = \widetilde{h}_d \label{lambp2}\\
 \max_{p\in [0,1]} \lambda_d(p) = h_d. \label{lambp3}
 \end{gather}
 \end{theo}
 \proof The proof of (\ref{lambp1}) is easy: pick any sequence
 $\m_q$ satisfying $\lim_{q \to \infty} \m_q =\infty$, and take $s_q =
 0$ for all $q$. Then conditions (\ref{HammCond}) hold for $p =
 0$, and consequently (\ref{lambpdef}) holds. But $\# W_0(\m_q,0)
 = 1$, since there is only one way to cover a box with monomers,
 and (\ref{lambp1}) follows.

 We prove (\ref{den}) and (\ref{lambp2}) together. Pick a sequence
 $(\m_q)_{q \in \N} \subseteq \N^d$ such that the $|\m_q|_{\pr}$
 are even and $\lim_{q\to\infty}\m_q =\infty$, and take $s_q =
 \frac{|\m_q|_{\pr}}{2}$. Then conditions (\ref{HammCond}) hold
 for $p = 1$, and consequently (\ref{lambpdef}) holds. But
 $W_0(\m_q,s_q) = \widetilde W_0(\m_q)$, and therefore
 \begin{equation}\label{star}
 \lim_{\m \to \infty, \frac{|\m|_{\pr}}{2}\in \N}
 \frac{\log \# \widetilde{W}_0(\m)}{|\m|_{\pr}} = \lambda_d(1).
 \end{equation}
 In view of (\ref{winb}) and (\ref{star}) we obtain
 \begin{multline}\label{star1}
 \lim_{\m \to \infty} \frac{\log \# \widetilde{W}(\m)}{|\m|_{\pr}} \geq
 \limsup_{\m \to \infty, \frac{|\m|_{\pr}}{2}\in \N}
 \frac{\log \# \widetilde{W}_{\perio}(\m)}{|\m|_{\pr}}
 \\
 \geq \liminf_{\m \to \infty, \frac{|\m|_{\pr}}{2}\in \N}
 \frac{\log \# \widetilde{W}_{\perio}(\m)}{|\m|_{\pr}}
 \geq \liminf_{\m \to \infty, \frac{|\m|_{\pr}}{2}\in \N}
 \frac{\log \# \widetilde{W}_0(\m)}{|\m|_{\pr}} = \lambda_d(1).
 \end{multline}
 For $\m \geq (2,\ldots,2) \in \N^d$ let
 $a(\m) := 2|\m|_{\pr}\sum_{i=1}^d\frac{1}{m_i}$ be the surface area of
 $\an{\m}$, and let
 \begin{gather*}
 w(\m) := \sum_{s \in [\frac{|\m|_{\pr}-a(\m)}{2},
 \frac{|\m|_{\pr}}{2}] \cap \Z } \#W_0(\m,s)\\
 \widetilde\omega(\m) := \max_{s\in [\frac{|\m|_{\pr}-a(\m)}{2},
 \frac{|\m|_{\pr}}{2}] \cap \Z} \#W_0(\m,s)\\
 \widetilde{s}(\m) := \argmax_{s\in [\frac{|\m|_{\pr}-a(\m)}{2},
 \frac{|\m|_{\pr}}{2}] \cap \Z} \#W_0(\m,s).
 \end{gather*}
 In words, $w(\m)$ is the sum of $\#W_0(\m,s)$ where $s$ ranges
 over those numbers of dimers that are sufficient to cover the interior of
 $\an{\m}$, i.e., the elements of $\an{\m}$ not on its boundary;
 $\widetilde{s}(\m)$ is the largest summand in that sum; and
 $\widetilde{s}(\m)$ is a number of dimers achieving the maximum.

 Clearly $\widetilde{\omega}(\m)\leq  w(\m) \leq
 \frac{a(\m) + 2}{2}\widetilde{\omega}(\m)$, and therefore
 \begin{equation}\label{star2}
 \limsup_{\m \to \infty} \frac{\log \widetilde\omega(\m)}{|\m|_{\pr}} =
 \limsup_{\m \to \infty} \frac{\log w(\m)}{|\m|_{\pr}}.
 \end{equation}
 By Proposition \ref{prop:multilimsup} there exists a sequence
 $(\n_q)_{q \in \N} \subseteq \N^d$ satisfying
 \begin{equation}\label{subsequence}
 \lim_{q \to \infty} \n_q = \infty, \quad
 \lim_{q \to \infty} \frac{\log \widetilde{\omega}(\n_q)}{|\n_q|_{\pr}} =
 \limsup_{\m \to \infty} \frac{\log \widetilde{\omega}(\m)}{|\m|_{\pr}}.
 \end{equation}
 Let $t_q := \widetilde{s}(\n_q)$ for each $q \in \N$,
 and so $\# W_0(\n_q,t_q) = \widetilde{\omega}(\n_q)$.
 Clearly $\lim_{q \to \infty} \frac{2 t_q}{|\n_q|_{\pr}}
 = 1$, and so conditions (\ref{HammCond}) hold for $\n_q,
 t_q$ with $p = 1$, and consequently (\ref{lambpdef})
 holds for them. Hence by (\ref{subsequence})
 \begin{equation}\label{star3}
 \limsup_{\m \to \infty} \frac{\log \widetilde{\omega}(\m)}{|\m|_{\pr}} =
 \lambda_d(1).
 \end{equation}
 Next we assert that $\#\widetilde{W}(\m) \leq  w(\m + 2\1)$.
 Indeed, each cover in $\widetilde{W}(\m)$ can be shifted by $\1$
 and extended by monomers to a tiling in $W_0(\m + 2\1,s)$ for one
 of the $s$ appearing in the sum $w(\m + 2\1)$. Therefore by
 (\ref{star2}) and (\ref{star3})
 \begin{multline}\label{star4}
 \lim_{\m \to \infty} \frac{\log \#\widetilde{W}(\m)}{|\m|_{\pr}} \leq
 \limsup_{\m \to \infty} \frac{\log w(\m + 2\1)}{|\m|_{\pr}} \\
 = \limsup_{\m \to \infty} \frac{\log w(\m)}{|\m|_{\pr}} =
 \limsup_{\m \to \infty} \frac{\log \widetilde{\omega}(\m)}{|\m|_{\pr}} =
 \lambda_d(1).
 \end{multline}
 Inequalities (\ref{star1}) and (\ref{star4}) combined, along with
 (\ref{winb}), complete the proof of (\ref{den}) and
 (\ref{lambp2}).

 We now prove (\ref{lambp3}). As $W_0(\m,s) \subseteq W_0(\m)$, it
 follows that $\lambda_d(p) \leq h_d$ for all $p \in [0,1]$. To complete the proof, we exhibit
 a $p^* \in [0,1]$ satisfying the reverse inequality.
 For each $\m \in \N^d$, let
 \begin{gather*}
 \omega(\m) := \max_{s \in [0,\frac{|\m|_{\pr}}{2}] \cap \Z}
 \#W_0(\m,s)
 \\
 s(\m) := \argmax _{s \in [0,\frac{|\m|_{\pr}}{2}] \cap \Z}
 \#W_0(\m,s)
 \\
 p(\m) := \frac{2s(\m)}{|\m|_{\pr}} \in [0,1],
 \end{gather*}
 so that $\omega(\m) = \#W_0(\m,s(\m))$.

 Observe that $\#W_0(\m) = \sum_{s \in [0,\frac{|\m|_{\pr}}{2}]
 \cap \Z} \#W_0(\m,s) \leq
 \frac{|\m|_{\pr} + 2}{2}\omega(\m)$, and therefore, by
 (\ref{dmen}),
 \begin{equation}\label{hdatmostliminf}
 h_d  \leq \liminf_{\m \to \infty} \frac{\log\omega(\m)}{|\m|_{\pr}}.
 \end{equation}
 From the bounded sequence $(p(q \1))_{q \in \N}$ choose a convergent subsequence
 $(p(q_k \1))_{k \in \N}$ and set $p^* := \lim_{k \to \infty} p(q_k \1) \in
 [0,1]$. Then conditions (\ref{HammCond}) hold for the sequences $q_k \1$
 and $s(q_k \1)$ with $p^*$, and therefore
 (\ref{lambpdef}) yields
 \begin{equation}\label{limislambdap*}
 \lim_{k \to \infty} \frac{\log\omega(q_k \1)}{q_k^d} = \lambda_d(p^*).
 \end{equation}
 By the definition of $\liminf$ we have $\liminf_{\m \to \infty}
 \frac{\log\omega(\m)}{|\m|_{\pr}} \leq \lim_{k \to \infty}
 \frac{\log\omega(q_k \1)}{q_k^d}$. Hence by
 (\ref{hdatmostliminf}) and (\ref{limislambdap*}) we obtain $h_d
 \leq  \lambda_d(p^*)$. \qed
 \begin{prop}\label{lowesth}  Let $d \in \N$.  Then for each $\m \in \N^d$
 \begin{gather}\label{lbw0}
 \frac{\log \# W(\m)}{|\m|_{\pr}} \geq h_d \geq \frac{\log \# W_0(\m)}{|\m|_{\pr}} \\
 \frac{\log \# \widetilde{W}(\m)}{|\m|_{\pr}} \geq \widetilde{h}_d \geq
 \frac{\log \# \widetilde{W}_0(\m)}{|\m|_{\pr}}.
 \end{gather}
 These upper and lower bounds converge to $h_d$ and
 $\widetilde{h}_d$, respectively, hence the latter are computable.
 \end{prop}
 \proof The upper bounds follow from the general theory of NNSOFT
 (\ref{ub1}), and their convergence from (\ref{hGamma}). For the
 lower bounds, let $k \in \N$ and consider the box $\an{k \m}$. It
 can be decomposed into $k^d$ shifted copies of $\an{\m}$. Hence
 \[\# W_0(k \m) \geq \# W_0(\m)^{k^d}, \quad
 \#\widetilde{W}_0(k \m) \geq \#\widetilde{W}_0(\m)^{k^d}.\]
  Sending $k$ to $\infty$ and using (\ref{dmen}) and (\ref{den}),
 we deduce the lower bounds as well as their convergence.  \qed

 We conclude this section by computing the various quantities in
 question for $d=1$ and illustrating Theorem \ref{maxlambp} for
 that case, where everything can be found explicitly. $\# W_0(m)$
 is the number of monomer-dimer tilings of $\an{m}$. Clearly it
 satisfies $\# W_0(1)=1$, $\# W_0(2)=2$ and $\# W_0(m) = \#
 W_0(m-1) + \# W_0(m-2)$ for $m \geq 3$. It follows that $\#
 W_0(m) = F_{m+1}$, where $F_m = \frac{1}{\sqrt{5}} \left(
 \frac{1+\sqrt{5}}{2} \right)^{m} - \frac{1}{\sqrt{5}} \left(
 \frac{1-\sqrt{5}}{2} \right)^{m}$ are the Fibonacci numbers. $\#
 W_{\perio}(m)$ is the number of monomer-dimer tilings of $T(m)$,
 and it satisfies $\# W_{\perio}(1) = 1$, $\# W_{\perio}(2) = 3$
 (one monomer tiling and two dimer tilings), and $\# W_{\perio}(m) =
 \# W_0(m) + \# W_0(m-2)$ for $m \geq 3$ (the second term counting
 the tilings with a dimer occupying $1$ and $m$). It follows that
 $\# W_{\perio}(m) = F_{m+1} + F_{m-1} = L_m$, where $L_m = \left(
 \frac{1+\sqrt{5}}{2} \right)^{m} + \left( \frac{1-\sqrt{5}}{2}
 \right)^{m}$ are the Lucas numbers. $\# W(m)$ is the number of
 monomer-dimer covers of $\an{m}$, where a dimer may protrude from
 $1$ to $0$ , or from $m$ to $m+1$. It satisfies $\# W(1) = 3$,
 $\# W(2) = 5$ and $\# W(m) = \# W_0(m) + 2 \# W_0(m-1)  + \#
 W_0(m-2)$ for $m \geq 3$ (the three terms representing covers
 with zero, one, or two protruding dimers, respectively). It
 follows that $\# W(m) = L_m + 2F_m = \left( 1 +
 \frac{2}{\sqrt{5}} \right) \left( \frac{1+\sqrt{5}}{2}
 \right)^{m} + \left( 1 - \frac{2}{\sqrt{5}} \right) \left(
 \frac{1-\sqrt{5}}{2} \right)^{m}$. From these values we see that
 $\frac{\log \# W(m)}{m}$,  $\frac{\log \# W_{\perio}(m)}{m}$ and
 $\frac{\log \# W_0(m)}{m}$ converge to $h_1 = \log
 \frac{1+\sqrt{5}}{2}$, in accordance with (\ref{dmen}).

 To determine $\lambda_1(p)$, it is enough to consider rational $p
 \in [0,1]$ by continuity, and then only $n \in \N$ such that $s =
 \frac{pn}{2} \in \N$, and send such $n$ to $\infty$,
 by(\ref{HammCond})--(\ref{lambpdef}). Then $\# W_0(s,n)$ is the
 number of linear arrangements of $s$ dimers and $n - 2s$
 monomers, which is equal to $\binom{n-s}{s}$. An application of
 Stirling's approximation then gives \cite{Ha1}
 \begin{multline*}
 \lambda_1(p)= \lim_{n \to \infty} \frac{1}{n} \log
 \binom{\left( 1  - \frac{p}{2} \right)n}{\frac{pn}{2}}
 \\
 = \left( 1 - \frac{p}{2} \right) \log \left(1 - \frac{p}{2} \right) - \frac{p}{2} \log
 \frac{p}{2} - (1 - p)\log (1 - p).
 \end{multline*}
 We see that $\lambda_1(0) = 0$ and $\lambda_1(1) = 0 =
 \widetilde{h}_1$ in accordance with (\ref{lambp1}) and
 (\ref{lambp2}). It is straightforward to verify that
 \[\max_{p \in [0,1]}\lambda_1(p) = \lambda_1\left( 1 - \frac{1}{\sqrt{5}} \right) = \log \frac{1+\sqrt{5}}{2} =h_1,\]
 in accordance with (\ref{lambp3}).

 \section{Lower Bounds for Monomer-Dimer Entropy with Dimer
 Density $p$}

 For an $m \times n$ matrix $A$, denote by $\pers A$ the sum of
 the permanents of all $s \times s$ submatrices of $A$. For a
 graph $G$, a \emph{matching} is a set of vertex-disjoint edges,
 and $W(G,s)$ denotes the set of all matchings of size $s$ in $G$,
 which can be regarded as covers of the vertex set $V(G)$ of $G$
 by $s$ dimers (edges) and $|V(G)| - 2s$ monomers (vertices). If
 $G$ is a bipartite graph with color classes $\an{m}$ and
 $\an{n}$, its \emph{incidence matrix} is the $m \times n$ $0$-$1$
 matrix $A = A(G)$ such that $a_{ij} = 1$ if and only if $\{i,j\}$ is
 an edge of $G$. In that case it is immediate that $\# W(G,s) =
 \pers A(G)$. A bipartite graph $G$ is said to be
 \emph{$r$-regular} if each vertex of $G$ has degree $r$,
 equivalently $A(G)$ has all row sums and column sums equal to
 $r$, so that $\frac{1}{r} A(G)$ is \emph{doubly-stochastic} (a
 nonnegative matrix with all row sums and column sums equal to
 $1$, necessarily a square matrix).

 \begin{theo}\label{lbmdcg}  Let $G$ be an $r$-regular bipartite
 graph with $n$ vertices in each color class.  Then
 \begin{equation}\label{lbmdcg1}
 \#W(G,s) \geq \binom{n}{s}^2 s!\,\left(\frac{r}{n}\right)^s.
 \end{equation}
 \end{theo}
 \proof A result of the first author \cite{Fr3} states that if $B$
 is a doubly-stochastic $n \times n$ matrix, then $\pers B \geq
 \pers J_n$, where $J_n$ is the $n \times n$ matrix with all
 entries equal to $\frac{1}{n}$. Since $\frac{1}{r} A(G)$ is
 doubly-stochastic, $\pers \frac{1}{r} A(G) =
 \frac{1}{r^s} \pers A(G)$ and $\pers J_n =
 \binom{n}{s}^2 \frac{s!}{n^s}$, the result follows. \qed

 The recent result of Schrijver \cite{Sch} improves this lower
 bound for the case $s=n$ if $r$ is constant and $n$ tends to
 infinity: under the assumptions of Theorem \ref{lbmdcg1}
 \begin{equation}\label{Schri}
 \#W(G,n) \geq  \left( \frac{(r-1)^{r-1}}{r^{r-2}} \right)^{n}.
 \end{equation}
 It would be of interest to similarly improve the lower bound of
 Theorem \ref{lbmdcg1} in the interesting range $n$ large and $s/n
 \geq r > 0$ (see below).

 In a recent paper \cite{Wan}, Wanless gives an alternative lower
 bound to (\ref{lbmdcg1}), namely
 $\#W(G,s) \geq \binom{n}{s} \left( \frac{(r-1)^{r-1}}{r^{r-2}}
 \right)^{s}$. It turns out that except for $\frac{s}{n}$ close to $1$,
 the bound (\ref{lbmdcg1}) is better.
 \begin{theo}\label{lpmdczd}  Let $ d \in \N$, $p \in [0,1]$ and recall the definition of
 $\lambda_d(p)$, the monomer-dimer entropy with dimer density $p$,
 given by (\ref{HammCond})--(\ref{lambpdef}).  Then
 \begin{equation}\label{lpmdczd1}
 \lambda_d(p) \geq \frac{1}{2} ( -p\log p - 2(1-p)\log(1-p) + p\log
 2d -p ).
 \end{equation}
 Furthermore, the dimer entropy $\widetilde{h}_d$ and monomer-dimer
 entropy $h_d$ satisfy
 \begin{equation}\label{lpmdczd2}
 \widetilde{h}_d=\lambda_d(1)  \geq  \frac{1}{2}((2d-1)\log(2d-1) -
 (2d-2)\log 2d),
 \end{equation}
 \begin{equation}\label{lpmdczd3}
 h_d  \geq  \frac{1}{2}(-p(d)\log p(d) - 2(1-p(d))\log(1-p(d)) +
 p(d)\log 2d -p(d)),
 \end{equation}
 where
 \begin{equation}\label{lpmdczd4}
 p(d)=\frac{4d+1-\sqrt{8d +1}}{4d}.
 \end{equation}
 \end{theo}
 \proof  Let $\m = (m_1,\ldots ,m_d)\in\N^d$ and assume that
 $m_1,\ldots,m_d$ are all even. Let $G$ be the adjacency graph of
 $T(\m)$. That is, the color classes of $G$ are the sets $\{ \i
 \in T(\m) : i_1 + \cdots + i_d\ \text{ even} \}$ and $\{ \j \in
 T(\m) : j_1 + \cdots + j_d\ \text{ odd} \}$, and $\{ \i,\j \}$ is an
 edge of $G$ if and only if $\i$ and $\j$ are neighbors on
 $T(\m)$, i.e., $\j = \i \pm \e_k$ for some $k \in \an{d}$, where
 the addition is the standard addition in the group
 $(\Z/m_1\Z)\times \cdots \times (\Z/m_d\Z)$.  Then $G$ is a
 $2d$-regular bipartite graph on $2n = |\m|_{\pr}$ vertices, and
 $W(G,s)$ is the set $W_{\perio}(\m,s)$ of monomer-dimer covers of
 $T(\m)$ having exactly $s$ dimers. Theorem \ref{lbmdcg} yields
 that $\#W_{\perio}(\m,s) \geq \binom{n}{s}^2 s! \left(
 \frac{2d}{n} \right)^s$. There is an injection $f$ from
 $W_{\perio}(\m,s)$ to $W_0(\m + \1,s)$, the set of monomer-dimer
 tilings of $\an{\m + \1}$ having exactly $s$ dimers. If $c \in
 W_{\perio}(\m,s)$, then $f(c)$ is obtained from $c$ by replacing
 each dimer in $c$ occupying the points $\i = (i_1,\ldots,i_d)$
 and $\j = \i + \e_k$ such that $i_k = m_k$ and $j_k = 0$ by a dimer
 occupying the points $\i$ and $(i_1,\ldots,i_{k-1},m_k +
 1,i_{k+1},\ldots,i_d)$. Therefore $\# W_0(\m + \1,s) \geq
 \binom{n}{s}^2 s! \left( \frac{2d}{n} \right)^s$. Let $(\m_q)_{q
 \in \N} \subseteq \N^d$ and $(s_q)_{q \in \N} \subseteq \N$ be
 sequences such that all the coordinates of each $\m_q$ are even,
 $\lim_{q \to \infty} \m_q = \infty$ and $\lim_{q \to \infty}
 \frac{2s_q}{|\m_q|_{\pr}} = p$. Set $n_q =
 \frac{|\m_q|_{\pr}}{2}$. Then conditions (\ref{HammCond}) hold,
 and consequently (\ref{lambpdef}) does. Therefore
 \begin{multline*}
 \lambda_d(p) = \lim_{q \to \infty} \frac{\log \# W_0(\m_q,s_q)}{|\m_q|}
 = \lim_{q \to \infty} \frac{\log \# W_0(\m_q + \1,s_q)}{|\m_q|}
 \\
 \geq \lim_{q \to \infty} \frac{\log \binom{n_q}{s_q}^2 s_q! \left( \frac{2d}{n_q}
 \right)^{s_q}}{2n_q}
 = \lim_{n \to \infty} \frac{1}{2n} \log \binom{n}{pn}^2 (pn)! \left( \frac{2d}{n}
 \right)^{pn}.
 \end{multline*}
 Manipulating the limit in the right-hand side of the inequality
 above and using the equality $\lim_{r \to \infty} \frac{1}{r}
 (\log r! - \log r^{r}) = -1$, we deduce the inequality
 (\ref{lpmdczd1}).

 Let $(\m_q)_{q \in \N}$ again satisfy the assumptions that all
 the coordinates of each $\m_q$ are even and $\lim_{q \to \infty}
 \m_q = \infty$, but this time set $s_q = n_q =
 \frac{|\m_q|_{\pr}}{2}$. Using the inequality (\ref{Schri}) for
 $\# W_{\perio}(\m_q,n_q)$ and (\ref{lambp2}), we
 deduce the inequality (\ref{lpmdczd2}).

 To prove (\ref{lpmdczd3}), we use (\ref{lambp3}). We easily
 verify that the right-hand side of (\ref{lpmdczd1}) is a strictly
 concave function of $p$ in $[0,1]$, and $p(d)$ given in
 (\ref{lpmdczd4}) is its unique critical point in that interval,
 hence its maximizing point there. \qed

 For $d = 2,3$, inequality (\ref{lpmdczd3}) yields
 \begin{align}
 h_2 &\geq 0.6358077435 \label{perlb2}\\
 h_3 &\geq 0.7652789557. \label{perlb3}
 \end{align}
 For $d=3$, inequality (\ref{lpmdczd2}) yields $\widetilde{h}_3
 \geq 0.440075842$, which is the best known lower bound.


 \section{Upper and Lower Bounds on $h_d$ and $\widetilde{h}_d$ Using Spectral Radii}

 For $d \in \N$, $K \subseteq \an{d}$ and $\m \in \N^d$, we denote
 by $\an{\m_K}$ the projection of $\an{\m}$ on the coordinates
 with indices in $K$. Let $W_{\perio,K}(\m)$, respectively
 $\widetilde{W}_{\perio,K}(\m)$, be the set of monomer-dimer
 covers, respectively dimer covers, of $T(\m_K) \times
 \an{\m_{\an{d} \setminus K}}$. Thus $W_{\perio,\an{d}}(\m) =
 W_{\perio}(\m)$ and $\widetilde{W}_{\perio,\an{d}}(\m) =
 \widetilde{W}_{\perio}(\m)$. Note that by the isotropy of our
 $\Gamma$, $\# W_{\perio,K}(\m)$ and $\# \widetilde{W}_{\perio,K}(\m)$
 are invariant under permutations of the components of $\m$ if $K$
 undergoes a corresponding change.

 In order to analyze $W_{\perio,\{d\}}(\m)$, we focus on the
 dimers in the cover lying along the direction $\e_d$. More
 precisely, with $\m' = (m_1,\ldots,m_{d-1})$, we consider
 $\an{\m'} \times T(m_d)$ as consisting of $m_d$ levels isomorphic
 to $\an{\m'}$. A subset $S$ of the points in level $q$ is covered
 by dimers joining levels $q-1$ and $q$ (with level $0$ understood
 as level $m_d$); a subset $T$ disjoint from $S$ is covered by
 dimers joining levels $q$ and $q+1$ (with level $m_d + 1$
 understood as level $1$); and the remainder $U$  of level $q$ is
 covered by monomers and dimers lying entirely within level $q$.
 We are interested in counting the coverings of $U$ subject to
 various restrictions. With that in mind, for $\m' \in \N^{d-1}$
 we define an undirected graph $G(\m')$ whose vertices are the
 subsets of $\an{\m'}$ in which subsets $S$ and $T$ are adjacent
 if and only if $S \cap T = \emptyset$. When $S \cap T =
 \emptyset$ we also define, using $U = \an{\m'} \setminus (S \cup
 T)$,
 \begin{align*}
   a_{ST} &= \text{ number of monomer-dimer tilings of } U
   \\
   b_{ST} &= \text{ number of monomer-dimer tilings of } U \text{ viewed as a subset of } T(\m')
   \\
   p_{ST} &=
   \begin{aligned}[t]
   &\text{ number of monomer-dimer covers of } U, \text{ viewed as a subset of }
   \\
   &T(m_1) \times \an{(m_2,\ldots,m_{d-1})},
   \text{ each monomer within } U, \text{ and each}
   \\
   &\text{ dimer meeting } U \text{ but not } S \cup T.
   \end{aligned}
   \\
   c_{ST} &=
   \begin{aligned}[t]
   &\text{ number of monomer-dimer covers of } U, \text{ each monomer
   within } U,
   \\
   &\text{ and each dimer meeting } U \text{ but not } S \cup T.
   \end{aligned}
 \end{align*}
 Thus in the tilings/covers counted by $a_{ST}$, $b_{ST}$,
 $p_{ST}$, $c_{ST}$, each monomer lies within $U$ and each dimer
 meets $U$ but not $S \cup T$. In $a_{ST}$, each dimer occupies
 two points of $U$ that are adjacent in $\an{\m'}$. In $b_{ST}$,
 each dimer occupies two points of $U$ that are adjacent in
 $T(\m')$, so is allowed to ``wrap around''. In $p_{ST}$, the
 dimers in the direction of $\e_1$ are allowed to ``wrap around''
 and the other dimers are allowed to ``protrude out'' of
 $\an{(m_2,\ldots,m_{d-1})}$. In $c_{ST}$, the dimers may
 ``protrude'' out of $\an{\m'}$. Therefore $a_{ST} \leq b_{ST}
 \leq p_{ST} \leq c_{ST}$. By definition, if $U = \emptyset$,
 then $a_{ST} = b_{ST} = p_{ST} = c_{ST} = 1$.
 Notice that when $d=2$, there is no distinction between $b_{ST}$
 and $p_{ST}$.

 We define the matrices $A(\m')$, $B(\m')$, $P(\m')$, $C(\m')$ with rows and
 columns indexed by subsets of $\an{\m'}$ as follows:

 \begin{align*}
   &\twocases{A(\m')_{ST}}{a_{ST}}{S \cap T = \emptyset}{0}{S \cap T \neq \emptyset}
   \\
   &\twocases{B(\m')_{ST}}{b_{ST}}{S \cap T = \emptyset}{0}{S \cap T \neq \emptyset}
   \\
   &\twocases{P(\m')_{ST}}{p_{ST}}{S \cap T = \emptyset}{0}{S \cap T \neq \emptyset}
   \\
   &\twocases{C(\m')_{ST}}{c_{ST}}{S \cap T = \emptyset}{0}{S \cap T \neq \emptyset.}
 \end{align*}
 Thus $A(\m')$, $B(\m')$, $P(\m')$, $C(\m')$ are symmetric matrices---here is
 the ``hidden symmetry'' referred to in Section 4---of integers
 satisfying $ 0 \leq A(\m') \leq B(\m') \leq P(\m') \leq C(\m')$ (where the
 inequalities indicate componentwise comparisons). We denote by
 $\alpha(\m'), \beta(\m'), \pi(\m'), \gamma(\m')$ their spectral radii,
 respectively, and consequently
 $\alpha(\m') \leq \beta(\m') \leq \pi(\m') \leq \gamma(\m')$.

 In an analogous way, we define $\widetilde{a}_{ST}$,
 $\widetilde{b}_{ST}$, $\widetilde{p}_{ST}$, $\widetilde{c}_{ST}$,
 where there are no monomers in the tilings and covers, the
 matrices $\widetilde{A}(\m')$, $\widetilde{B}(\m')$,
 $\widetilde{P}(\m')$, $\widetilde{C}(\m')$ and their spectral
 radii $\widetilde{\alpha}(\m'), \widetilde{\beta}(\m'),
 \widetilde{\pi}(\m'), \widetilde{\gamma}(\m')$.

 Each of these eight symmetric matrices can be considered as the
 adjacency matrix of an undirected multigraph, where the
 multiplicity of an edge is the corresponding matrix entry. This
 multigraph is a weighted version of $G(\m')$. If the multigraph
 is bipartite, we say that the matrix is \emph{bipartite}; if the
 multigraph is connected, we say that the matrix is
 \emph{irreducible}; if the multigraph is disconnected, we say
 that the matrix is a \emph{direct sum}; if the multigraph is
 connected and the greatest common divisor of the lengths of all
 its closed walks is 1, we say that the matrix is \emph{primitive},
 equivalently for sufficiently high powers of the matrix, all
 entries are strictly positive.

 \begin{prop}\label{symprop}
 Let $2 \leq  d\in \N$ and $\m = (\m',m_d) \in \N^d$. Then
 \begin{enumerate}
  \item[(a)] $\tr A(\m')^{m_d}$ is the number of monomer-dimer
  tilings of $\an{\m'} \times T(m_d)$ and $\tr \widetilde{A}(\m')^{m_d}$
  is the number of dimer tilings of $\an{\m'} \times T(m_d)$;
  \item[(b)] $\tr B(\m')^{m_d} =
  \# W_{\perio}(\m)$
  and $\tr \widetilde{B}(\m')^{m_d} =
  \# \widetilde{W}_{\perio}(\m)$;
  \item[(c)] $\tr P(\m')^{m_d} =
  \# W_{\perio,\{1,d\}}(\m)$
  and $\tr \widetilde{P}(\m')^{m_d} =
  \# \widetilde{W}_{\perio,\{1,d\}}(\m)$;
  \item [(d)] $\tr C(\m')^{m_d} = \#
  W_{\perio,\{d\}}(\m)$ and $\tr \widetilde{C}(\m')^{m_d} = \#
  \widetilde{W}_{\perio,\{d\}}(\m)$;
  \item[(e)] for $m_d \geq 2$, if column vector $\x = (x_S)_{S \subseteq \an{\m'}}$
  is given by $x_S = b_{S\emptyset}$, then
  $\x^{\trans} B(\m')^{m_d - 2} \x = \# W_{\perio,\an{d-1}}(\m)$,
  if vector $\y$ is given by $y_S = c_{S\emptyset}$, then
  $\y^{\trans} C(\m')^{m_d - 2} \y = \# W(\m)$,
  and if $\z = (z_S)_{S \subseteq \an{\m'}}$
  is given by $z_S = p_{S\emptyset}$, then
  $\z^{\trans} P(\m')^{m_d - 2} \z = \# W_{\perio,\{1\}}(\m)$;
  if column vector $\widetilde{\x} = (\widetilde{x}_S)_{S \subseteq \an{\m'}}$
  is given by $\widetilde{x}_S = \widetilde{b}_{S\emptyset}$, then
  $\widetilde{\x}^{\trans} \widetilde{B}(\m')^{m_d - 2} \widetilde{\x}
  = \# \widetilde{W}_{\perio,\an{d-1}}(\m)$,
  if $\widetilde{\y}$ is given by $\widetilde{y}_S = \widetilde{c}_{S\emptyset}$, then
  $\widetilde{\y}^{\trans} \widetilde{C}(\m')^{m_d - 2} \widetilde{\y}
  = \# \widetilde{W}(\m)$,
  and if vector $\widetilde{\z}$ is given by $\widetilde{z}_S = \widetilde{p}_{S\emptyset}$,
  then $\widetilde{\z}^{\trans} \widetilde{P}(\m')^{m_d - 2} \widetilde{\z}
  = \# \widetilde{W}_{\perio,\{1\}}(\m)$;
  \item[(f)] the matrices $A(\m')$, $B(\m')$, $P(\m')$, $C(\m')$ are primitive;
  \item[(g)] if $|\m'|_{\pr}$ is odd, then $\widetilde{A}(\m')$,
  $\widetilde{B}(\m')$ are bipartite, otherwise they are direct sums.
 \end{enumerate}
 \end{prop}
 \proof We begin with proving the first part of (b), its second
 part and (a), (c), (d) and (e) being similar. Assume first that $m_d =
 1$, and let $\phi \in W_{\perio}(\m)$. Since $\phi$ can be
 extended periodically in the direction of $\e_d$ with period 1,
 it can be viewed as an element of $W_{\perio}(\m')$. Therefore
 $\# W_{\perio}(\m) = \# W_{\perio}(\m')$. We have $\tr B(\m') =
 \sum_{S \subseteq \an{\m'}} b_{SS}$. Only the term $S =
 \emptyset$ contributes to the sum, and for this term we have $U =
 \an{\m'}$ and $b_{\emptyset\emptyset} = \# W_{\perio}(\m')$.
 Hence $\tr B(\m') = \# W_{\perio}(\m')$. Now assume that $m_d >
 1$, and consider any closed path $S_1,S_2,\ldots,S_{m_d},S_1$ of
 length $m_d$ in $G(\m')$. For each $\p' \in S_q$ place a dimer
 occupying the points $(\p',q)$ and $(\p',q+1)$ (with $m_d + 1$
 wrapping around to $1$). We want to extend these dimers to a
 monomer-dimer tiling of $T(\m') \times T(m_d) = T(\m)$, i.e., to
 a member of $W_{\perio}(\m)$, by monomers and by dimers not in
 the direction of $\e_d$, i.e., lying within the levels
 $1,\ldots,m_d$. The number of choices of such monomers and dimers
 to fill the remainder of level $q$ is given by $b_{S_{q-1}S_q}$,
 and so the number of extensions to a member of $W_{\perio}(\m)$
 is $b_{S_{1}S_{2}} b_{S_{2}S_{3}} \cdots b_{S_{m_d - 1}S_{m_d}}
 b_{S_{m_d}S_{1}}$. Conversely, each member of $W_{\perio}(\m)$ is
 obtained in this way. Hence $\# W_{\perio}(\m)$ is the sum of all
 the products of the above form, namely $\tr B(\m')^{m_d}$.

 To prove (f), we note that $A(\m')$ is irreducible, since
 whenever $S \cap T = \emptyset$, $U$ can be tiled by monomers and
 therefore each subset of $\an{\m'}$ is adjacent to $\emptyset$ in
 the graph of $A(\m')$. Furthermore, $A(\m')$ is primitive since
 the graph has a cycle of length 1 from $\emptyset$ to
 $\emptyset$. Since $A(\m') \leq B(\m') \leq P(\m') \leq C(\m')$,
 $B(\m')$, $P(\m')$ and $C(\m')$ are also primitive.

 To prove (g), let $\cE, \cO$ denote the subsets of $\an{\m'}$
 with even and odd cardinality, respectively. If
 $\widetilde{b}_{ST} > 0$, then $U$ can be tiled by dimers and so
 $\# U$ must be even. Therefore if $|\m'|_{\pr}$ is odd, members
 of $\cE$ are adjacent only to members of $\cO$ in the graph of
 $\widetilde{B}(\m')$, and so that graph is bipartite; if
 $|\m'|_{\pr}$ is even, then members of $\cE$ are adjacent only to
 themselves, and the graph is disconnected. The same conclusions
 hold for $\widetilde{A}(\m')$ since $\widetilde{A}(\m') \leq
 \widetilde{B}(\m')$.
 \qed
 \begin{lemma}\label{albetgain}
 Let $2 \leq  d\in \N$ and $\m' \in \N^{d-1}$.  Then
 \begin{align}
  &\lim_{m_d \to \infty} \frac{\log \# W_0(\m',m_d)}{m_d} =
  \log \alpha(\m')\label{albetgaina}
  \\
  &\lim_{m_d \to \infty}
  \frac{\log \# W_{\perio,\an{d-1}}(\m',m_d)}{m_d}=
  \log \beta(\m')\label{albetgainb}
  \\
  &\lim_{m_d \to \infty}
  \frac{\log \# W_{\perio,\{1\}}(\m',m_d)}{m_d}=
  \log \pi(\m')\label{albetgainbprime}
  \\
  &\lim_{m_d \to \infty} \frac{\log \# W(\m',m_d)}{m_d}=
  \log \gamma(\m')\label{albetgainc}
  \\
  &\lim_{m_d \to \infty} \frac{\log \#\widetilde
  W_0(\m',m_d)}{m_d} \leq \log \widetilde\alpha(\m')\label{albetgaind}
  \\
  &\lim_{m_d \to \infty}\frac{\log \#
  \widetilde{W}_{\perio,\an{d-1}}(\m',m_d)}{m_d}=
  \log \widetilde\beta(\m')\label{albetgaine}
  \\
  &\lim_{m_d \to \infty}\frac{\log \#
  \widetilde{W}_{\perio,\{1\}}(\m',m_d)}{m_d}=
  \log \widetilde\pi(\m')\label{albetgaineprime}
  \\
  &\lim_{m_d \to \infty}
  \frac{\log \# \widetilde{W}(\m',m_d)}{m_d}=
  \log \widetilde\gamma(\m').\label{albetgainf}
  \end{align}
  \end{lemma}
  \proof From Part (a) of Proposition \ref{symprop} we obtain
  $\#\widetilde{W}_0(\m',m_d)
  \leq \tr \widetilde{A}(\m')^{m_d}$, and therefore
  \begin{equation}\label{limsuplelog}
  \limsup_{m_d \to \infty} \frac{\log \# \widetilde{W}_0(\m',m_d)}{m_d}
  \leq
  \limsup_{m_d \to \infty} \frac{\log \tr \widetilde{A}(\m')^{m_d}}{m_d} = \log \widetilde{\alpha}(\m').
  \end{equation}
  The equality in (\ref{limsuplelog}) follows from a
  characterization of $\rho(M)$ for a square matrix $M \geq 0$,
  namely $\rho(M) = \limsup_{n \to \infty} (\tr
  M^n)^{\frac{1}{n}}$ (see for example Proposition 10.3 of
  \cite{Fr2}). Since $-\log \# \widetilde{W}_0(\m',m_d)$ is
  subadditive in $m_d$, the first $\limsup$ in (\ref{limsuplelog})
  can be replaced by a $\lim$, which proves (\ref{albetgaind}).
  Similar considerations prove $\lim_{m_d \to \infty} \frac{\log
  \# W_0(\m',m_d)}{m_d} \leq \log \alpha(\m')$. In order to prove
  the reverse inequality and thus (\ref{albetgaina}), observe that
  each monomer-dimer tiling of $\an{\m'} \times T(m_d)$ extends to a
  monomer-dimer tiling in $W_0(\m',m_d + 1)$ (replace each dimer
  occupying $(\m',1)$ and $(\m',m_d)$ by a monomer occupying
  $(\m',1)$ and a dimer occupying $(\m',m_d)$ and $(\m',m_d + 1)$,
  and tile the rest with monomers). Hence $\# W_0(\m',m_d + 1)
  \geq \tr A(\m')^{m_d}$ by Part (a) of Proposition \ref{symprop}.
  Therefore, since $-\log \# W_0(\m',m_d)$ is subadditive in $m_d$
  and thus the limits below exist, we obtain
  \begin{multline*}
  \lim_{m_d \to \infty} \frac{\log \# W_0(\m',m_d)}{m_d}
  = \lim_{m_d \to \infty} \frac{\log \# W_0(\m',m_d + 1)}{m_d}
  \\
  \geq \limsup_{m_d \to \infty} \frac{\log \tr A(\m')^{m_d}}{m_d}
  = \log \alpha(\m').
  \end{multline*}

  To prove (\ref{albetgainb}), (\ref{albetgainbprime}) and
  (\ref{albetgainc}), we use another characterization of the
  spectral radius. A \emph{vector norm} is a mapping $\| \cdot \|
  : M_n(\C) \rightarrow \R_{+}$ taking complex matrices of order
  $n$ to nonnegative reals such that $\| M \| = 0$ only if $M =
  0$, $\| zM \| = |z| \| M \|$ for all $z \in \C$, and $\| M + N
  \| \leq \| M \| + \| N \|$. If $c_{ij} > 0$ for all $i,j \in
  \an{n}$, then $\| M \| = \sum_{ij} c_{ij} |m_{ij}|$ is a vector
  norm. Proposition 10.1 of \cite{Fr2} states that if $\| \cdot
  \|$ is a vector norm, then $\rho(M) = \lim_{k \to \infty} \| M^k
  \|^{\frac{1}{k}}$. In particular, if $M \geq 0$ and $\v$ is a
  column vector with positive entries, then $\rho(M) = \lim_{k \to
  \infty} (\v^{\trans} M^k \v)^{\frac{1}{k}}$. Applying this to $M
  = B(\m'), P(\m'), C(\m')$ and using Part (e) of Proposition
  \ref{symprop} with $\v = \x,\z,\y$ defined there proves
  (\ref{albetgainb}), \ref{albetgainbprime}), (\ref{albetgainc}).

  The proof of (\ref{albetgaine}) is a little more complicated
  because the vector $\widetilde{\x}$ in Part (e) of Proposition
  \ref{symprop} is not strictly positive. Therefore we introduce
  the vector $\widetilde{\w}$ with entries $\widetilde{w}_S =
  \max(1, \widetilde{x}_S)$. Then, by Part (e) of Proposition
  \ref{symprop}, we have $\# \widetilde{W}_{\perio,\an{d-1}}(\m) =
  \widetilde{\x}^{\trans} \widetilde{B}(\m')^{m_d - 2}
  \widetilde{\x} \leq \widetilde{\w}^{\trans}
  \widetilde{B}(\m')^{m_d - 2} \widetilde{\w}$. Therefore we
  obtain
  \[\lim_{m_d \to \infty} \frac{\log \#  \widetilde{W}_{\perio,\an{d-1}}(\m)}{m_d}
  \leq \lim_{m_d \to \infty} \frac{\log \widetilde{\w}^{\trans}
  \widetilde{B}(\m')^{m_d - 2} \widetilde{\w}}{m_d} = \log
  \widetilde{\beta}(\m')\] (the first $\lim$ above exists since
  $\log \# \widetilde{W}_{\perio,\an{d-1}}(\m)$ is subadditive in
  $m_d$). On the other hand $\#
  \widetilde{W}_{\perio,\an{d-1}}(\m) \geq \#
  \widetilde{W}_{\perio}(\m) = \tr \widetilde{B}(\m')^{m_d}$ by
  Part (b) of Proposition \ref{symprop}. Therefore
  \[\lim_{m_d \to \infty} \frac{\log \# \widetilde{W}_{\perio,\an{d-1}}(\m)}{m_d}
  \geq \limsup_{m_d \to \infty} \frac{\log \tr \widetilde{B}(\m')^{m_d}}{m_d}
  = \log \widetilde{\beta}(\m').\]
  This proves (\ref{albetgaine}). To prove (\ref{albetgainf}), we
  show analogously that
  \[\lim_{m_d \to \infty} \frac{\log \#  \widetilde{W}(\m)}{m_d}
  \leq \log \widetilde{\gamma}(\m'),\]
  and on the other hand, by Part (d) of Proposition \ref{symprop},
  \begin{multline*}
  \lim_{m_d \to \infty} \frac{\log \# \widetilde{W}(\m)}{m_d}
  \geq \limsup_{m_d \to \infty} \frac{\log \# \widetilde{W}_{\perio,\{d\}}(\m)}{m_d}
  \\
  = \limsup_{m_d \to \infty} \frac{\log \tr \widetilde{C}(\m')^{m_d}}{m_d}
  = \log \widetilde{\gamma}(\m').
  \end{multline*}
  The proof of (\ref{albetgaineprime}) is similar.
  \qed

  \begin{prop}\label{stanlub}
   Let $2 \leq  d\in \N$ and $\m'\in \N^{d-1}$.  Then
   \begin{alignat}{2}
   \frac{\log \gamma(\m')}{|\m'|_{pr}} &\geq\; &h_d
   &\geq  \frac{\log \alpha(\m')}{|\m'|_{pr}} \label{stanlub1}
   \\
   \frac{\log \widetilde\gamma(\m')}{|\m'|_{pr}}&\geq\; &\widetilde{h}_d
   &\geq \frac{\log \widetilde\alpha(\m')}{|\m'|_{pr}}. \label{stanlub2}
   \end{alignat}
   \end{prop}
   \proof  The upper bounds follow from the general upper bounds
   in Proposition \ref{lowesth} along with (\ref{albetgainc}),
   (\ref{albetgainf}). The lower bound in (\ref{stanlub1}) follows
   similarly from the general lower bound in Proposition
   \ref{lowesth} along with (\ref{albetgaina}). However, since
   (\ref{albetgaind}) only gives a lower bound for $\log
   \alpha(\m')$, we use a separate argument for the lower bound
   in (\ref{stanlub2}) as follows. For $q \in \N$, $\#
   \widetilde{W}(q\m',m_d)$ is not smaller than the the number of
   dimer tilings of $\an{q\m'} \times T(m_d)$, which in turn is
   not smaller than the number of dimer tilings of $\an{\m'} \times
   T(m_d)$ raised to the $q^{d-1}$ power. Hence by Part (a) of
   Proposition \ref{symprop} we have
   \[\# \widetilde{W}(q\m',m_d) \geq \left( \tr \widetilde{A}^{m_d} \right)^{q^{d-1}},\]
   and so
   \[\frac{\log \# \widetilde{W}(q\m',m_d)}{|(q\m',m_d)|_{\pr}} \geq
   \frac{\log \tr \widetilde{A}^{m_d}}{|\m'|_{\pr} m_d}.\]
   Therefore
   \[\widetilde{h}_d = \lim_{q,m_d \to \infty} \frac{\log \# \widetilde{W}(q\m',m_d)}{|(q\m',m_d)|_{\pr}}
   \geq \frac{1}{|\m'|_{\pr}} \limsup_{q,m_d \to \infty} \frac{\log \tr \widetilde{A}^{m_d}}{m_d}
   = \frac{\log \widetilde\alpha(\m')}{|\m'|_{pr}}.\]
  \qed

  Now we introduce the following notation. For $\m \in \N^d$ and
  $k \in \an{d}$, $\m^{\sim k} := (m_1,\ldots,m_{k-1},
  m_{k+1},\ldots,m_d)\in \N^{d-1}$. As special cases we have the
  previous notation $\m' = \m^{\sim d}$ and $\m^- = \m^{\sim 1}$.
  \begin{prop}\label{perest}  Let $\m \in \N^d$,
  and assume that $m_d$ is even.  Then each
  $k \in \an{d-1}$ satisfies
  \begin{align}
  \frac{\log \beta(\m^{\sim d})}{|\m|_{\pr}} &\leq
  \frac{\log 2}{m_k} + \frac{\log \beta(\m^{\sim k})}{|\m^{\sim k}|_{pr}} \label{perest1}
  \\
  \frac{\log \widetilde{\beta}(\m^{\sim d})}{|\m|_{\pr}} &\leq
  \frac{\log 2}{m_k} + \frac{\log \widetilde{\beta}(\m^{\sim k})}{|\m^{\sim k}|_{pr}}. \label{perest2}
  \end{align}
  \end{prop}
  \proof We have
  \begin{equation*}
  \begin{split}
  \beta(\m^{\sim d})^{m_d} \leq \tr B(\m^{\sim d})^{m_d} = \# W_{\perio}(\m)
  = \tr B(\m^{\sim k})^{m_k}
  \\
  \leq 2^{|\m^{\sim k}|_{\pr}} \beta(\m^{\sim k})^{m_k},
  \end{split}
  \end{equation*}
  where the first inequality follows since $\beta(\m^{\sim d})$ is
  one of the eigenvalues of $B(\m^{\sim d})$, which are all real,
  and $m_d$ is even, the next equality from Part (b) of
  Proposition \ref{symprop}, the next equality from the same and
  the fact that $\# W_{\perio}(\m)$ is invariant under coordinate
  permutations in $\m$, and the last inequality from the fact that
  $B(\m^{\sim k})$ has $2^{|\m^{\sim k}|_{\pr}}$ eigenvalues, all
  real, whose absolute values are at most $\beta(\m^{\sim k})$.
  Taking logarithms and dividing by $|\m|_{pr}$, we deduce
  (\ref{perest1}). The inequality (\ref{perest2}) is obtained in a
  similar way.
  \qed
  We define
  \begin{align*}
  \overline{h}_{d-1}(m_1) &:= \lim_{\m^{-} \to \infty}
  \frac{\log \# W_{\perio,\{1\}}(m_1,\m^{-})}
  {|\m^{-}|_{\pr}},\; m_1 \in \N;
  \quad \overline{h}_{d-1}(0):= \log 2
  \\
  \breve{h}_{d-1}(m_1) &:= \lim_{\m^{-} \to \infty}
  \frac{\log \# \widetilde{W}_{\perio,\{1\}}(m_1,\m^{-})}
  {|\m^{-}|_{\pr}},\; m_1 \in \N;
  \quad \breve{h}_{d-1}(0):= \log 2.
  \end{align*}
  Notice that for $m_1 \in \N$, $\overline{h}_{d-1}(m_1)$ is
  the same as $\overline{h}(m_1,\Gamma)$ defined in
  (\ref{defbarh}) when $\Gamma$ is the $d$-digraph encoding the
  monomer-dimer covers. For this reason the limit
  $\overline{h}_{d-1}(m_1)$ exists, and similarly for
  $\breve{h}_{d-1}(m_1)$.
  The following theorem is an analog of Theorem \ref{perulb} and Theorem \ref{ubspecrt}.
  \begin{theo}\label{main}  Let $2 \leq  d \in \N,\;p,r\in \N,\;q\in \Z_+$.  Then
  \begin{alignat}{2}
  \frac{\overline{h}_{d-1}(2r)}{2r} &\geq\;  &h_d &\geq
  \frac{\overline{h}_{d-1}(p+2q)-\overline{h}_{d-1}(2q)}{p}
  \label{imlb}\\
  \frac{\breve{h}_{d-1}(2r)}{2r} &\geq\; &\widetilde{h}_d &\geq
  \frac{\breve{h}_{d-1}(p+2q)-\breve{h}_{d-1}(2q)}{p}. \label{imlh}
  \end{alignat}
  Let $\m'=(m_1,\ldots,m_{d-1}) \in \N^{d-1}$ and assume that $m_1,\ldots,m_{d-1}$
  are even.  Then
  \begin{align}
  h_d &\leq \frac{\beta(\m')}{|\m'|_{pr}} \label{imub1}
  \\
  \widetilde{h}_d &\leq \frac{\widetilde\beta(\m')}{|\m'|_{pr}}.  \label{imub2}
  \end{align}
  \end{theo}

  \proof
  We have
  \begin{multline*}
  \widetilde{h}_d = \lim_{\stack{\m',m_d \to \infty}{\frac{m_d}{2} \in \N}}
  \frac{\log \# \widetilde{W}_0(\m',m_d)}{|\m'|_{\pr}m_d}
  \leq \liminf_{\m' \to \infty}\frac{\log \widetilde{\alpha}(\m')}{|\m'|_{\pr}}
  \leq \limsup_{\m' \to \infty}\frac{\log \widetilde{\gamma}(\m')}{|\m'|_{\pr}}
  \\
  = \limsup_{\m',m_d \to \infty} \frac{\log \# W(\m',m_d)}{|\m'|_{\pr}m_d}
  = \widetilde{h}_d,
  \end{multline*}
  where the first equality follows from (\ref{den}), the next
  inequality from (\ref{albetgaind}), the next one from
  $\widetilde{\alpha}(\m') \leq \widetilde{\gamma}(\m')$, the next
  equality from (\ref{albetgainf}), and the last equality again
  from (\ref{den}). From this and
  $\widetilde{\alpha}(\m') \leq \widetilde{\beta}(\m') \leq \widetilde{\gamma}(\m')$
  we obtain
  \begin{equation}\label{stcharhd}
  \widetilde{h}_d = \lim_{\m' \to \infty} \frac{\log \widetilde{\alpha}(\m')}{|\m'|_{\pr}}
      = \lim_{\m' \to \infty} \frac{\log \widetilde{\beta}(\m')}{|\m'|_{\pr}}
      = \lim_{\m' \to \infty} \frac{\log \widetilde{\gamma}(\m')}{|\m'|_{\pr}}.
  \end{equation}
  Similarly (and more simply)
  \begin{equation}\label{stcharhmd}
  h_d = \lim_{\m' \to \infty} \frac{\log \alpha(\m')}{|\m'|_{\pr}}
      = \lim_{\m' \to \infty} \frac{\log \beta(\m')}{|\m'|_{\pr}}
      = \lim_{\m' \to \infty} \frac{\log \gamma(\m')}{|\m'|_{\pr}}.
  \end{equation}

  First we prove (\ref{imub1}).
  Let $\m' = (m_1,\ldots,m_{d-1}) \in \N$, $m_1,\ldots,m_{d-1}$
  even, and let $\p = (p_1,\ldots,p_{d-1}) \in \N^{d-1}$ be arbitrary.  Set
  \begin{multline*}
  \m_1 = (p_1,\ldots,p_{d-1},m_1),\quad \m_2 = (p_2,\ldots,p_{d},m_1,m_2),\quad \ldots,
  \\
  \m_{d-1} = (p_{d},m_1,\ldots,m_{d-1}).
  \end{multline*}
  Then, using (\ref{perest1}) with $k = 1$ $d-1$ times,
  we obtain
  \begin{multline*}
  \frac{\log \beta(\p)}{|\p|_{pr}} \leq \frac{\log 2}{p_1} + \frac{\log \beta(\m_1^-)
  }{|\m_1^-|_{pr}} \leq \frac{\log 2}{p_1} + \frac{\log 2}{p_2} +
  \frac{\log\beta(\m_2^-)}{|\m_2^-|_{pr}} \leq \cdots
  \\
  \leq \sum_{j=1}^{d-1}\frac{\log 2}{p_j} + \frac{\log \beta(\m')}
  {|\m'|_{pr}}.
  \end{multline*}
  Letting $\p\to\infty$ and using (\ref{stcharhmd}) for the
  left-hand side, we deduce (\ref{imub1}). Similar arguments apply
  to deduce (\ref{imub2}).

  We now demonstrate the lower bound in (\ref{imlb}). Let $\m^-
  \in \N^{d-1}$, $p \in \N$, $q\in \Z_+$. Assume first that $q \in
  \N$. Since $\gamma(\m^-) = \rho(C(\m^-))$ and $C(\m^-)$ is
  symmetric, it follows as in the arguments for
  (\ref{rhogeqtraceratio}) that
  \begin{equation}\label{lbhd}
  \gamma(\m^-)^p \geq \frac{\tr C(\m^-)^{p+2q}}{\tr C(\m^-)^{2q}} =
  \frac{\# W_{\perio,\{1\}}(p+2q,\m^-)}{\# W_{\perio,\{1\}}(2q,\m^-)}.
  \end{equation}
  Taking logarithms, dividing by $|\m^-|_{pr}$, letting $\m^- \to
  \infty$, and using (\ref{stcharhmd}) and the definition of
  $\overline{h}_{d-1}(m_1)$, we deduce the lower bound in (\ref{imlb}) for the
  case $q \in \N$. If $q = 0$, we have to replace the denominators
  in (\ref{lbhd}) by $\tr I = 2^{|\m^-|_{\pr}}$, and the lower
  bound in (\ref{imlb}) is verified because $h_{d-1}(0)$ was
  defined to be $\log 2$. The lower bound in (\ref{imlh}) is proved
  similarly.

  We now prove the upper bound of (\ref{imlb}).  For each $\m' \in
  \N^{d-1}$ we have
  \[\gamma(\m')^{2r} \leq \tr C(\m')^{2r} = \# W_{\perio,\{d\}}(\m',2r )
  = \#W_{\perio,\{1\}}(2r,\m'),\]
  where the inequality above is
  true because the eigenvalues of the symmetric matrix $C(\m')$
  are real and $\gamma(\m')$ is one of them, the first equality
  follows from Part (d) of Proposition \ref{symprop}, and the last
  equality from the invariance under coordinate permutations.
  Therefore
  \[\frac{\log \gamma(\m')}{|\m'|_{pr}} \leq
  \frac{\log \#W_{\perio,\{1\}}(2r,\m')}{2r|\m'|_{pr}},\] and
  letting $\m' \to \infty$, we deduce the upper bound of
  (\ref{imlb}) by (\ref{stcharhmd}) and the definition of
  $\overline{h}_{d-1}(m_1)$. Similarly we deduce the upper bound
  of (\ref{imlh}). \qed

  The following theorem supplies practical upper and lower bounds
  on $2$- and $3$-dimensional monomer-dimer and dimer entropies.
  \begin{theo}\label{dmcor}  Let $p,r,t,u,v \in \N$ and $q,s\in \Z_+$.  Then
  \begin{gather*}
  \frac{\log \beta(2r)}{2r} \geq h_2
  \geq \frac{\log \beta(p+2q) -\log \beta(2q)}{p}, \quad \beta(0) = 2
  \\
  \frac{\log\tilde\beta(2r)}{2r} \geq \widetilde{h}_2 \geq
  \frac{\log \widetilde{\beta}(p+2q) - \log \widetilde{\beta}(2q)}{p},
  \quad \widetilde{\beta}(0) = 2
  \\
  \frac{\log \beta(2r,2t)}{4rt} \geq h_3
  \geq \frac{\log\beta(p+2q,u+2s)
   - \log\beta(p+2q,2s)}{up} - \frac{\log \beta(2q,2v)}{2vp}
  \\
  \frac{\log\tilde\beta(2r,2t)}{4rt}
  \geq \widetilde{h}_3
  \geq \frac{\log \widetilde{\beta}(p+2q,u+2s) - \log \widetilde{\beta}(p+2q,2s)}{up} -
  \frac{\log \widetilde{\beta}(2q,2v)}{2vp}
  \\
  \beta(n,0) = \beta(0,n) = \widetilde{\beta}(n,0) = \widetilde{\beta}(0,n) = 2^n, \quad n \in \N.
  \end{gather*}
  \end{theo}
 \proof  The upper bounds in the above inequalities are the
 inequalities (\ref{imub1}) and (\ref{imub2}).  We now show the
 lower bounds. Equations (\ref{albetgainb}) and (\ref{albetgaine})
 for $d = 2$ yield
 \begin{equation}\label{hbarlogbeta}
  \overline{h}_1(m_1) = \log \beta(m_1), \quad
  \breve{h}_1(m_1) = \log \widetilde{\beta}(m_1), \quad m_1 \in \N.
 \end{equation}
  Hence the lower bounds on $h_2, \widetilde{h}_2$ follow immediately from
  the lower bounds in
 (\ref{imlb}), (\ref{imlh}), equation(\ref{hbarlogbeta}) and the equalities
 $\overline{h}_1(0) = \breve{h}_1(0) = \log 2$.

 In order to establish the lower bounds on on $h_3,
 \widetilde{h}_3$, we first establish lower and upper bounds on
 $\overline{h}_2(m_1)$ and $\breve{h}_2(m_1)$ in terms of
 $\beta(\cdot,\cdot)$ and $\widetilde\beta(\cdot,\cdot)$.
 The definition of $\overline{h}_2(m_1)$ and $\breve{h}_2(m_1)$
 and equations (\ref{albetgainbprime}) and (\ref{albetgaineprime}) for $d = 3$ yield
 \begin{equation}\label{prophd}
 \overline{h}_2(m_1) = \lim_{m_2 \to \infty} \frac{\log \pi(\m')}{m_2},
 \quad \breve{h}_2(m_1 )= \lim_{m_2 \to \infty}\frac{\log \widetilde{\pi}(\m')}{m_2},
 \quad m_1 \in \N,
 \end{equation}
 where $\m' = (m_1,m_2)$.
 Since $P(\m')$ is a nonnegative symmetric matrix with spectral
 radius $\pi(\m')$, it follows as in (\ref{rhogeqtraceratio}) and using
 Part (c) of Proposition \ref{symprop} that
 \[\pi(\m')^u\ge \frac{\tr P(\m')^{u+2s}}{\tr P(\m')^{2s}}=\frac{\#W_{\perio,\{1,3\}}(\m',u+2s)}
 {\#W_{\perio,\{1,3\}}(\m',2s)}.\]
 Here $u\in \N$ and $s\in \Z_+$.
 When $s = 0$, $\tr P(\m')^{2s} =  2^{|\m'|_{\pr}}$, and so this
 is the value we use for
 $\#W_{per,\{1,3\}}(\m',0)$. Take logarithms of this inequality, divide by
 $m_2$ and send $m_2$ to $\infty$.  Using (\ref{prophd}) and (\ref{albetgainb}) for $d=3$, we deduce that
 \begin{equation}\label{prophd1}
 \overline{h}_2(m_1) \geq \frac{\log \beta(m_1,u+2s) - \log \beta(m_1,2s)}{u},
 \quad m_1 \in \N,
 \end{equation}
 where $\beta(m_1,0) := 2^{m_1}$. Similarly
 \begin{equation}\label{prophd2}
 \breve{h}_2(m_1) \geq \frac{\log \widetilde{\beta}(m_1,u+2s) - \log \widetilde{\beta}(m_1,2s)}{u},
 \quad m_1 \in \N,
 \end{equation}
 where $\widetilde{\beta}(m_1,0) : = 2^{m_1}$.
 For $v\in \N$ we have the inequality
 $\pi(\m')^{2v} \leq \tr P(\m')^{2v} = \#W_{\perio,\{1,3\}} (\m',2v)$.
 Take logarithms of this inequality,
 divide by $2vm_2$ and send $m_2$ to $\infty$. Using (\ref{prophd}) and
 (\ref{albetgainb}) for $d=3$, we deduce that for $m_1 \in \N$
 \begin{equation}\label{prophd3}
 \overline{h}_2(m_1) \leq \frac{\log \beta(m_1,2v)}{2v}.
 \end{equation}
 Inequality (\ref{prophd3}) also holds for $m_1 = 0$ since by
 definition $\overline{h}(0) = \log 2$ and $\beta(0,2v) = 2^{2v}$.
 Similarly, for $m_1 \in \Z_+$
 \begin{equation}\label{prophd4}
 \breve{h}_2(m_1) \le \frac{\log \widetilde{\beta}(m_1,2v)}{2v}.
 \end{equation}

 Now we can substitute the bounds (\ref{prophd1}) and
 (\ref{prophd3}) in the lower bound of (\ref{imlb}) as appropriate
 from the signs in the numerator, and obtain the lower bound on
 $h_3$ as stated in the theorem, and similarly for
 $\widetilde{h}_3$. \qed

 \section{Using Automorphism Subgroups to Reduce Computations}

 The matrix $B(\m')$ has order $2^n$, where $n =
 |\m'|_{\pr}$, and so has $4^n$ entries. Since its $(S,T)$
 entries are positive precisely when $S \cap T = \emptyset$, its
 number of positive entries is $\sum \binom{n}{i} 2^{n-i}
 = 3^n$. Hence it is sparse. However, already for $\m' = (4,4)$
 it has $4.3 \cdot 10^7$ nonzero entries, and the computation
 of its spectral radius is infeasible for standard PC. Nevertheless,
 this computation can be reduced to computing the spectral radii
 of a suitable nonnegative matrix whose order is the number of
 orbits of the action of an automorphism subgroup of $B(\m')$.
 This usage of automorphisms is also used in \cite{Ciu} and \cite
 {Lun}.

 Recall that given an $N \times N$ complex-valued matrix $A =
 (a_{ij})_1^N$, its automorphism group is the subgroup of the
 symmetric group $S_N$ on $\an{N}$ defined by
 \begin{equation}\label{autgrp}
  \Aut(A) := \{ \pi \in S_N \;:\; a_{\pi(i)\pi(j)}= a_{ij} \text{ for all } i,j \in \an{N}\}.
 \end{equation}
 Let $\cG$ be a subgroup of $\Aut(A)$. The action of $\cG$
 partitions $\an{N}$ into minimal invariant subsets called orbits.
 We denote by $\cO := \an{N}/\cG$ the orbit space (set of orbits),
 and by Greek letters $\alpha,\beta,\ldots$ its members. We have
 \begin{equation}\label{invariance}
  \sum_{j \in \beta} a_{ij} = \sum_{j \in \beta} a_{\pi(i)\pi(j)}
  = \sum_{k \in \beta} a_{\pi(i)k}, \qquad \alpha,\beta \in \cO,\;\; i
  \in \alpha,\;\; \pi \in S_N,
 \end{equation}
 which means that for given $\alpha,\beta \in \cO$, the sum $\Sigma_{j \in \beta}
 a_{ij}$ is the same for all $i \in \alpha$.
 Let $M=\#\cO$, and define the $M \times M$ matrix
 $\widehat{A}=(\widehat{a}_{\alpha\beta})_{\alpha,\beta \in \cO}$
 by
 \begin{equation}\label{defahat}
 \widehat{a}_{\alpha\beta}= \sum_{j \in \beta} a_{ij}, \quad i\in\alpha.
 \end{equation}
 This is a valid definition by (\ref{invariance}).
 The following proposition is known, and we prove it for completeness.
 \begin{prop}\label{invspec}
 Let $A=(a_{ij})_1^N$ be a complex-valued matrix. Let $\cG$ be a
 subgroup of $\Aut A$, $\cO$ its orbit space, and $M = \#\cO$. Let
 $\widehat{A}$ be the induced $M \times M$ complex-valued matrix
 given by (\ref{defahat}). Then the spectrum (set of eigenvalues)
 of $\widehat{A}$, $\spec(\widehat{A})$, is a subset of
 $\spec(A)$, and in particular $\rho(\widehat{A}) \leq \rho(A)$. If
 $A$ is a real-valued nonnegative matrix, then $
 \rho(\widehat{A}) = \rho(A)$. If $A$ is real and symmetric, then
 $\widehat{A}$ is symmetric with respect to an appropriate inner
 product on $\R^M$, and in particular $\spec(\widehat{A})$ is real
 and $\widehat{A}$ is diagonalizable.
 \end{prop}
 \proof  Let $\Pi_N$ be the group of $N \times N$ permutation matrices.  Let
 $\iota: S_N \to \Pi_N$ be the standard representation of $S_N$. That is
 $\iota(\pi)(x_i)_{i \in \an{N}} =  (x_{\pi(i)})_{i\in \an{N}}$.  Let
 \begin{multline*}
 \cX := \{ \x \in \C^N : \iota(\pi)(\x) = \x
 \text{ for all } \pi \in \cG \}
 \\
 = \{(x_i)_{i \in \an{N}} \in \C^N: x_{\pi(i)} = x_i
 \text{ for all } i \in \an{N}, \pi \in \cG \}
 \end{multline*}
 be the subspace of vectors that are constant on each orbit of
 $\cG$. Then $\cX \subseteq \C^N$ is the largest subspace of
 $\C^N$ on which $\iota(\cG)$ acts trivially (as the identity
 operator).  Clearly, $\cX$ is isomorphic to $\C^M$. Indeed, each
 $\x = (x_i)\in \cX$ induces a unique vector
 $\widehat{\x}:=(\widehat{x}_{\alpha})_{\alpha \in \cO} \in \C^M$,
 where $\widehat{x}_{\alpha} = x_i$ for any $i\in\alpha$.
 Conversely, each $\y \in \C^M$ induces a unique
 $\x \in \cX$ such that $\y = \widehat{\x}$. Next, we observe that
 $\cX$ is an invariant subspace of $A$. Indeed, for each $\x =
 (x_i) \in \cX$ and $\pi \in \cG$ we have for all $i \in \an{N}$
 \[
 (A\x)_i = \sum_{j=1}^N a_{ij}x_j
 = \sum_{j=1}^N a_{\pi(i)\pi(j)} x_{\pi(j)}
 = \sum_{k=1}^N a_{\pi(i)k}x_k
 = (A\x)_{\pi(i)},
 \]
 which means that $A\x \in \cX$.
 Moreover, if  $\x \in \cX$ and $\widehat{\x} =
 (\widehat{x}_{\alpha}) \in \C^M$ is defined as above, then for
 any $i \in \alpha$ we have $(A\x)_i=\sum_{\beta\in \cO} \widehat
 a_{\alpha\beta} \widehat{x}_{\beta}$, and consequently
 $\widehat{A\x}=\widehat{A}\widehat{\x}$.
 This means that the action of $A|_{\cX}$ is isomorphic the the
 action of $\widehat A$ on $\C^M$. In particular,
 \[\spec(\widehat{A}) = \spec(A|_{\cX}) \subseteq \spec(A),\]
 and therefore
 \[\rho(\widehat{A})  \leq  \rho(A).\]
 Assume now that $A$ is nonnegative. Then by the Perron-Frobenius
 theorem $\rho(A) \in \spec(A)$, and $A$ has an eigenvector $\x$
 belonging to $\rho(A)$. Since each $\pi \in \Aut(A)$ satisfies $A
 \iota(\pi) = \iota(\pi) A$, it follows that $\iota(\pi) \x$ is
 also an eigenvector of $A$ belonging to $\rho(A)$.  Hence
 $\sum_{\pi \in \Aut(A)} \iota(\pi)\x \in \cX$ is an eigenvector
 of $A$ belonging to $\rho(A)$. Therefore $\rho(A) \in
 \spec(A|_{\cX}) = \spec(\widehat{A})$. It follows that
 $\rho(\widehat A) = \rho(A)$.

 Finally assume that $A$ is a real symmetric matrix.  That is
 $(A\x,\y)=(\x,A\y)$, where $(\x,\y) = \y^{\trans}\x$ is the standard inner
 product in $\R^N$. For each $\alpha \in \cO$, let $w_{\alpha}$ be
 the cardinality of the orbit $\alpha$.  In $\R^M$ we define the inner
 product
 \begin{equation}\label{newinp}
 \an{\widehat{\x},\widehat{\y}} := \sum_{\alpha \in \cO}
 w_{\alpha} \widehat{x}_{\alpha} \widehat{y}_{\alpha}.
 \end{equation}
 Then all $\x,\y \in \cX$ satisfy $(\x,\y) =
 \an{\widehat{\x},\widehat{\y}}$.  Hence
 $\an{\widehat{A}\widehat{\x},\widehat{\y}} =
 \an{\widehat{\x},\widehat{A}\widehat{\y}}$, i.e., $\widehat{A}$
 is symmetric (self adjoint) with respect to the inner product
 $\an{\cdot,\cdot}$ in $\R^M$. In particular, $\widehat{A}$ has
 real eigenvalues and is similar to a diagonal matrix.  \qed

 We shall now briefly mention the power method for computing
 $\rho(A)$ where $A$ is a nonnegative symmetric matrix of order
 $N$, and a variant of it that works on $\widehat{A}$ of
 order $M$, which we used in our computations.
 \begin{prop}\label{power}
 Let $A$ be a nonnegative symmetric matrix of order $N$.
 Choose a scalar $r > 0$ and a positive vector
 $\x_0 = (x_{0,1},\ldots,x_{0,N})^{\trans}$.
 For each $m \in \N$, let
 \begin{gather*}
  \x_m = (x_{m,1},\ldots,x_{m,N})^{\trans} := (A + rI)\x_{m-1}\\
  l_m := \min_i \frac{x_{m,i}}{x_{m-1,i}}\\
  u_m := \max_i \frac{x_{m,i}}{x_{m-1,i}}\\
  r_m := \frac{(\x_m,\x_{m-1})}{(\x_{m-1},\x_{m-1})}.
 \end{gather*}
 Then $l_m$ is nondecreasing and $u_m$ is nonincreasing in $m$,
 \begin{gather*}
 l_m \leq r_m \leq \rho(A) + r \leq u_m, \qquad m \in \N \\
 \lim_{m \to \infty} l_m = \lim_{m \to \infty} u_m = \rho(A) + r,
 \end{gather*}
 and $\x_m/\sqrt{(\x_m,\x_m)}$ converges to an eigenvector
 of $A$ belonging to $\rho(A)$.

 Furthermore, with the notation of Proposition \ref{invspec}, if
 we choose the vector $\x_0$ to be in $\cX$, i.e., if
 $\x_0$ is constant on each orbit of $\cG$, then each $m
 \in \N$ the vector $\x_m$  is also in $\cX$ (so $\widehat{\x}_m$ is defined),
 \begin{gather*}
  \widehat{\x}_m = (\widehat{A} + r
  \widehat{I})\widehat{\x}_{m-1}\\
  l_m = \min_{\alpha \in \cO} \frac{\widehat{x}_{m,\alpha}}{\widehat{x}_{m-1,\alpha}}\\
  u_m = \max_{\alpha \in \cO} \frac{\widehat{x}_{m,\alpha}}{\widehat{x}_{m-1,\alpha}}\\
  r_m =
  \frac{\an{\widehat{\x}_m,\widehat{\x}_{m-1}}}{\an{\widehat{\x}_{m-1},\widehat{\x}_{m-1}}},
 \end{gather*}
  and $\widehat{\x}_m/\sqrt{(\widehat{\x}_m,\widehat{\x}_m)}$ converges to an eigenvector
  of $\widehat{A}$ belonging to $\rho(\widehat{A}) = \rho(A)$.
 \end{prop}

 For $\m' \in \N^{d-1}$, let $G_T(\m')$ be the adjacency graph of
 the elements of the torus $T(\m')$. The automorphisms of
 $G_T(\m')$ act as automorphisms of the symmetric nonnegative
 matrices $B(\m')$ and $\widetilde{B}(\m')$. In view of
 Proposition \ref{power}, in order to compute the spectral radii
 $\beta(\m')$ and $\widetilde{\beta}(\m')$, it is advantageous to
 use large automorphism subgroups of $G_T(\m')$. The rigid motions
 of the box $\an{\m'}$ and of the torus $T(\m')$ are
 automorphisms of $G_T(\m')$.

 The rigid motions of $\an{\m'}$ contain the reflections across
 the hyperplanes $x_k = \frac{m_k + 1}{2}$, $k \in \an{d-1}$,
 which commute with each other, and the allowable transpositions
 exchanging $x_i$ and $x_j$ in case $m_i = m_j$. Thus if $\m' = m
 \1$, $m \geq 2$, then the group of rigid motions of the cube
 $\an{\m'}$ contains a subgroup of order $2^{d-1}(d-1)!$. For $d =
 2,3$, which is our main focus in this paper, the reflections and
 allowable transpositions generate all the rigid motions of
 $\an{\m'}$.

 The rigid motions of $T(\m')$ contain, in addition to the rigid
 motions of $\an{\m'}$, the unit translations $\x \mapsto \x +
 \e_k$, $k \in \an{d-1}$. The unit translations generate the group
 of translations, an Abelian group isomorphic to $(\Z/m_1\Z)
 \times \cdots \times (\Z/m_{d-1}\Z)$ of order $|\m'|_{\pr}$. We
 call the group generated by the reflections, the allowable
 transpositions and the unit translations \emph{the group of rigid
 motions of $T(\m')$}. Note that for $T(2)$ the reflection
 coincides with the unit translation, and similarly for $T(\m')$,
 if $m_k = 2$ then the reflection across $x_k = \frac{3}{2}$
 coincides with the unit translation $\x \mapsto \x + \e_k$. We
 are aware of additional automorphisms  of $G_T(\m')$ if at least
 two components of $\m'$ are equal to 4: observe that $G_T(4)$ is
 isomorphic to $G_T(2,2)$, since both are $4$-cycles. Therefore
 $G_T(4,4)$ is isomorphic to $G_T(2,2,2,2)$, and its automorphism
 group has order at least $2^4 \cdot 4! = 384$, whereas the group
 of rigid motions of $T(4,4)$ has order $2^2 \cdot 2 \cdot 4^2 =
 128$. Similar results hold for $d > 3$.

 The following proposition is straightforward:
 \begin{prop}\label{trangrp}
 Let $\Gamma_1, \ldots, \Gamma_d \subseteq \an{n} \times \an{n}$
 and $\Gamma = (\Gamma_1, \ldots, \Gamma_d)$. Let $\m = (m_1,
 \ldots, m_d) \in \N^d$ and $\m'=(m_1, \ldots, m_{d-1})$, and
 consider the transfer digraph $\Theta_d(\m')$ between members of
 $W_{\perio}(\m')$ with respect to $\Gamma_d$. Then the group of
 translations of $T(\m')$ acts a subgroup of automorphisms of
 $\Theta_d(\m')$. If for some $k \in \an{d-1}$ $\Gamma_k$ is
 symmetric, then the reflection across the hyperplane $x_k =
 \frac{m_k+1}{2}$ acts as an automorphism of $\Theta_d(\m')$.  If
 for some $p,q \in \an{d-1}$ $m_p = m_q$ and $\Gamma_p =
 \Gamma_q$, then the transposition exchanging $x_p$ and $x_q$ acts
 as an automorphism of $\Theta_d(\m')$.
 \end{prop}
 \begin{corol}\label{autsub}
 Let $\Gamma_1, \ldots, \Gamma_d \subseteq \an{n} \times \an{n}$,
 $\Gamma = (\Gamma_1, \ldots, \Gamma_d)$, and assume that
 $\Gamma_1, \ldots, \Gamma_{d-1}$ are symmetric. Let $\m = (m_1,
 \ldots, m_d) \in \N^d$ and $\m'=(m_1, \ldots, m_{d-1})$, and
 assume that for all $p,q \in \an{d-1}$, $\Gamma_p =
 \Gamma_q$ if $m_p = m_q$. Then the automorphism subgroup of
 $G_T(\m')$ described above acts as an automorphism subgroup of the
 transfer digraph $\Theta_d(\m')$.
 \end{corol}

 As an example, consider the upper and lower bounds given by
 (\ref{ubspecr2}). The parameter $\theta_2(m)$ appearing there is
 the spectral radius of the matrix $B(m)$ defined in Section 6,
 which has an automorphism subgroup of order $2m$, isomorphic to
 the group of rigid motions of $T(m)$, if $m > 2$. $B(15)$ is
 $2^{15} \times 2^{15}$, but as we shall see, $\widehat{B}(15)$ is
 $ 1224 \times 1224 $, which makes the computation of its spectral
 radius feasible on a regular desktop computer.

 These observations are our main keys in finding good upper and
 lower bounds for $h_2$ and $h_3$. We point out that \cite{Ciu}
 was the first work that used these automorphisms of
 $\widetilde{B}(\m')$ to help obtain a good upper bound for
 $\widetilde{h}_3$, which was later improved in \cite{Lun} by
 similar methods.

 \section{Numerical Results for Monomer-Dimer Entropy in Two and Three Dimensions}

 Our results are based on Theorem \ref{dmcor}, and we compute the
 spectral radii appearing there using Propositions \ref{invspec}
 and \ref{power}, and the automorphism subgroups described in
 Section 7. We first consider the two-dimensional monomer-dimer
 entropy. Recall that $\beta(m_1)$ is the spectral radius of
 $B(m_1)$.  Table \ref{nresh2} lists $\log \beta (m_1)$, $\frac{\log
 \beta(m_1)}{m_1}$, and the number $\#\cO(m_1)$ of orbits of the
 torus $T(m_1)$ under the action of the group of rigid motions of
 $T(m_1)$.
 \begin{table}
 \[
  \begin{array}{cccc}
   m_1&\# \cO(m_1)&\log\beta(m_1)&\frac{\log\beta(m_1)}{m_1}\\
   \hline
   4 & 6 & 2.6532941163 & .66332352908 \\
   5 & 8 & 3.3135066910 & .66270133821 \\
   6 & 13 & 3.9769139475 & .66281899125 \\
   7 & 18 & 4.6395628723 & .66279469604 \\
   8 & 30 & 5.3023993987 & .66279992338 \\
   9 & 46 & 5.9651887945 & .66279875494 \\
   10 & 78 & 6.6279902386 & .66279902386 \\
   11 & 126 & 7.2907885674 & .66279896067 \\
   12 & 224 & 7.9535877093 & .66279897578 \\
   13 & 380 & 8.6163866375 & .66279897212 \\
   14 & 687 & 9.2791856222 & .66279897301 \\
   15 & 1224 & 9.9419845918 & .66279897279 \\
   16 & 2250 & 10.60478356551861 & .662798972844913 \\
   17 & 4112 & \in (11.26758254,11.26758315) & \in (.6627989729,.6627990088)
  \end{array}
 \]
 \caption{Spectral radii for $h_2$}\label{nresh2}
 \end{table}
 The computation of $\log\beta(17)$ was interrupted, and the table
 indicates the best interval in which we can locate it. We notice
 that the sequence $\frac{\log\beta(2r)}{2r}$ is decreasing for $r
 = 2, \ldots, 8$. Hence $h_2 \leq \frac{\log\beta(16)}{16} =
 .662798972844913$ is the best upper bound for $h_2$ from our
 data. The best lower bound for $h_2$ from our data is
 $h_2 \geq \frac{\log\beta(17) - \log\beta(16)}{1} \geq .66279897$. This
 improves the lower bound (\ref{perlb2}) from permanents by more
 than 4\%. Hence we obtain the value
 \begin{equation}\label{cnvh2}
 h_2=.66279897,
 \end{equation}
 correct to 8 decimal digits.  We also notice that the sequence
 $\frac{\log\beta(2j+1)}{2j+1}$ is increasing for $j=2, \ldots,
 8$. Suppose that this sequence were increasing for all values of
 $j$. Since $\lim_{j \to \infty} \frac{\log\beta(2j+1)}{2j+1} =
 h_2$ by (\ref{stcharhmd}), it would follow that $h_2 \geq
 \frac{\log\beta(17)}{17} \geq .6627989729$. The last digit of
 this bound is too high, as seen by comparison with our best upper
 bound, probably caused by roundoff errors in the interrupted
 computation, but enables us to state that the above hypothesis
 would gives the value $h_2 = .6627989728$ correct to $10$ digits,
 consistent with the one found by Baxter \cite{Bax1} (his value of
 $h_2$ is accurate to 8 digits, as can be seen by evaluating $\log
 \frac{\kappa}{s}$ for $s=1$ in his Table II and varying the last
 digit of the tabulated $\frac{\kappa}{s}$). Since the lower bound
 (\ref{perlb2}) for $h_2$ is quite close to the correct value of
 $h_2$, it is reasonable to assume that the value $p^*$, for which
 $\lambda_2(p^*) = h_2$, is fairly close to $p(2) =
 \frac{9-\sqrt{17}}{8} \sim 0.6096118$ (according to \cite{Bax1}, $p^*=0.63812311$.).

 As a check, Table \ref{nresth2} gives $\widetilde{\beta}(m_1)$, the
 spectral radius of $\widetilde{B}(m_1)$, yielding lower and upper bounds
 for the known entropy $\widetilde{h}_2 = 0.29156090\ldots$.
 \begin{table}
 \[
  \begin{array}{cccc}
   m_1&\# \cO(m_1)&\log\widetilde{\beta}(m_1)&\frac{\log\widetilde\beta(m_1)}{m_1}\\
   \hline
   4   &                6 &                                  1.316957897   &  .3292 \\
   5   &                8 &                                  1.404661127   &  .2809 \\
   6   &               13 &                                  1.843797237   &  .3073 \\
   7   &               18 &                                  2.003260294   &  .2862 \\
   8   &               30 &                                  2.400842203   &  .3001 \\
   9   &               46 &                                  2.594837310   &  .2883 \\
   10   &               78 &                                  2.969359257   &  .2969 \\
   11   &              126 &                                  3.183303939   &  .2894 \\
   12   &              224 &                                  3.543130579   &  .2953 \\
   13   &              380 &                                  3.770113562   &  .2900 \\
   14   &              687 &                                  4.119721251   &  .2943 \\
   15   &             1224 &                                  4.355934472   &  .2904
  \end{array}
 \]
 \caption{Spectral radii for $\widetilde{h}_2$}\label{nresth2}
 \end{table}
 Again, the sequence $\frac{\log\widetilde{\beta}(2r)}{2r}$ decreases
 for $r = 2, \ldots, 7$ and the sequence
 $\frac{\log\widetilde{\beta}(2j+1)}{2j+1}$ increases for $j=2,
 \ldots, 7 $. Thus the best upper bound on $\widetilde{h}_2$ from
 our data is $\frac{\log\widetilde{\beta}(14)}{14} = .2943$, which is
 larger by $0.9\%$ than the true value. The best lower bound is
 $\frac{\log\widetilde{\beta}(14) - \log\widetilde{\beta}(12)}{2} = 0.2883$,
 which is smaller by $1.1\%$ than the true value. We notice that
 $\frac{\log\widetilde{\beta}(15)}{15} = .2905 < \widetilde{h}_2$,
 consistent with the assumed fact that
 $\frac{\log\widetilde{\beta}(2j+1)}{2j+1}$ increases for all $j$.

 We now consider the three-dimensional monomer-dimer entropy
 $h_3$. Recall that $\beta(m_1,m_2) = \beta(m_2,m_1)$ is the
 spectral radius of $B(m_1,m_2)$.  Table \ref{nresh} gives $\log\beta (m_1,m_2)$,
 $\frac{\log\beta(m_1,m_2)}{m_1m_2}$, and the number $\#\cO(m_1,m_2)$
 of orbits of the torus $T(m_1,m_2)$ under the action of the group
 of rigid motions of $T(m_1,m_2)$.  (In the case
 $(m_1,m_2)=(4,4)$, we recall that the group of rigid motion of
 $T(2,2)$ has order $128$, and it turns out to have 805 orbits.
 We also did the computations with the larger automorphism
 subgroup of $G_T(4,4)$ of order $384$ discussed in Section 7,
 which turns out to have 402 orbits. Both computations gave the
 same value of $\beta(4,4)$.)
 \begin{table}
 \[
  \begin{array}{cccc}
   (m_1,m_2)&\# \cO(m_1,m_2)&\log\beta(m_1,m_2)&\frac{\log\beta(m_1,m_2)}{m_1m_2}\\
   \hline
   (2,2) & 6 & 3.224405658 &  0.8061014145 \\
   (3,2) & 13  & 4.768958913 & 0.7948264855 \\
   (4,2) & 34 &   6.367778959 &  0.7959723699 \\
   (5,2) & 78 & 7.958105292  & 0.7958105292 \\
   (6,2) & 237 & 9.550024542 & 0.7958353785 \\
   (7,2) &  687  & 11.14163679 & 0.7958311993 \\
   (8,2) &  2299  & 12.73331093 & 0.7958319331 \\
   (3,3) & 25  &  7.057039652 &   0.7841155169 \\
   (4,3) & 158 &   9.421594940 &  0.7851329117 \\
   (5,3) & 708 & 11.77517604 &  0.7850117360 \\
   (4,4) & 805 & 12.57923752 & 0.7862023450 \\
  \end{array}
 \]
 \caption{Spectral radii for $h_3$}\label{nresh}
 \end{table}
 Recall that $h_3 \leq \frac{\log\beta(2r,2t)}{4rt}$, and hence
 the best upper bound for $h_3$ from our data is
 $\frac{\log\beta(4,4)}{16} = 0.7862023450$. The best lower bound
 is $\frac{\log\beta(3,5) - \log\beta(3,4)}{1 \cdot
 1}-\frac{\log\beta(2,8)}{8 \cdot 1}=.761917234$. It turns out
 that the permanent lower bound (\ref{perlb3}) is better: $h_3
 \geq .7652789557$. Of course, had we computed $\beta(m_1,m_2)$
 for larger $m_1$ and $m_2$, we would eventually improve the
 permanent lower bound. Thus, the best estimates we have are
 \begin{equation}\label{cnvh3}
 .7652789557 \leq h_3 \leq .7862023450.
 \end{equation}
 Table \ref{nresth} lists $\widetilde{\beta}(m_1,m_2)$, the spectral
 radius of $\widetilde{B}(m_1,m_2)$, which give bounds for
 $\widetilde{h}_3$. The entry $(m_1,m_2) = (6,4)$ is taken from
 \cite{Lun}, which took advantage of the fact that the matrix of
 order $184854$ is a direct sum of 3 matrices.
 \begin{table}
 \[
  \begin{array}{cccc}
   (m_1,m_2)&\# \cO(m_1,m_2)&\log\widetilde{\beta}(m_1,m_2)&\frac{\log\widetilde{\beta}(m_1,m_2)}{m_1m_2}\\
   \hline
   (2,2) & 6  & 2.292431670 &  0.5731079175 \\
   (3,2) &13  & 3.068671222 &  0.5114452037 \\
   (4,2) &34  & 4.151763891 &  0.5189704864 \\
   (5,2) &78  & 5.119835223 &  0.5119835223 \\
   (6,2) &237 & 6.161467494 &  0.5134556245 \\
   (7,2) &687   & 7.168058989 &  0.5120042135 \\
   (3,3) &25  & 3.938705096 &  0.4376338996 \\
   (4,3) &158 & 5.365527945 &  0.4471273287 \\
   (5,3) &708 & 6.635849120 & 0.4423899413 \\
   (4,4) &805 & 7.409698288 &  0.4631061430\\
   (6,3) &4236   & 7.97716207        &  0.443175671\\
   (6,4) &184854   & 10.98112634 &  0.4575469308\\
  \end{array}
 \]
 \caption{Spectral radii for $\widetilde{h}_3$}\label{nresth}
 \end{table}
 The best upper bound for $\widetilde{h}_3$ is
 $\frac{\log\widetilde{\beta}(6,4)}{6 \cdot 4} = 0.4575469308$,
 which was reported in \cite{Lun}.  The best lower bound from the
 data is given by $\frac{\log\widetilde{\beta}(4,6) -
 \log\widetilde{\beta}(4,4)}{2 \cdot 2} -
 \frac{\log\widetilde{\beta}(2,6)}{6 \cdot 2} = .3794013885$,
 which is a weak lower bound.  The best lower bound for
 $\widetilde{h}_3$ is given by (\ref{lpmdczd2}): $\widetilde{h}_3
 \geq 0.4400758$.

 We now compare our results for $h_2$ with the results of
 \cite{HaM}.  On page 342, Hammersley and Menon tabulate estimates
 of $\lambda_2(p)$ computed by the Monte Carlo method in
 increments of $0.05$ for $ 0 \leq p \leq 1$.  The maximal value
 in their table is $.6676$ for $p = 0.65$. They state ``There are
 reasons for believing that this Monte Carlo estimate has a small
 negative bias, probably $1\%$ or $2\%$ too low".  However, since
 $\lambda_2(p) \leq h_2 = .66279897$, the Monte Carlo estimate for
 $\lambda_2(0.65)$ is at least $0.7\%$ higher than the true
 value.

\begin{figure}
  \centering
  \includegraphics[clip=true,angle=270,scale=0.55]{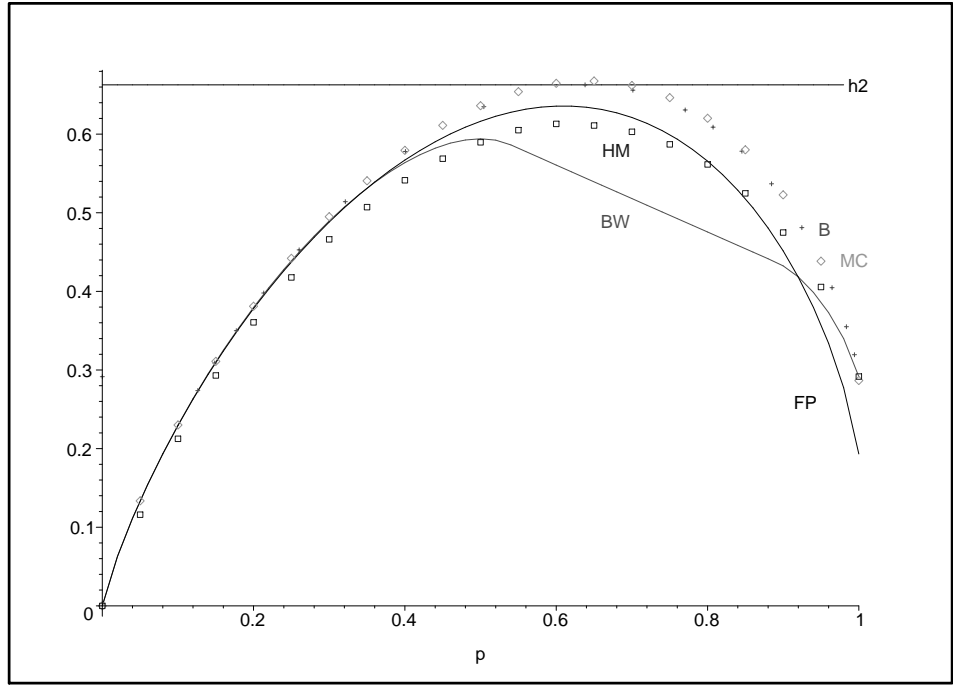}
  \caption{Lower bounds and estimates for $\lambda_2(p)$. HM is
  the lower bound of \cite{HaM}, BW is the lower bound of
  \cite{BoW}, FP is the lower bound of Theorem
  \protect\ref{lpmdczd}, MC is the Monte Carlo estimate of
  \cite{HaM}, B is the estimate from \cite{Bax1}, and h2 is the
  true value of $h_2 = \max \lambda_2(p)$.} \label{figbounds2}
 \end{figure}

\begin{figure}
  \centering
  \includegraphics[clip=true,angle=270,scale=0.55]{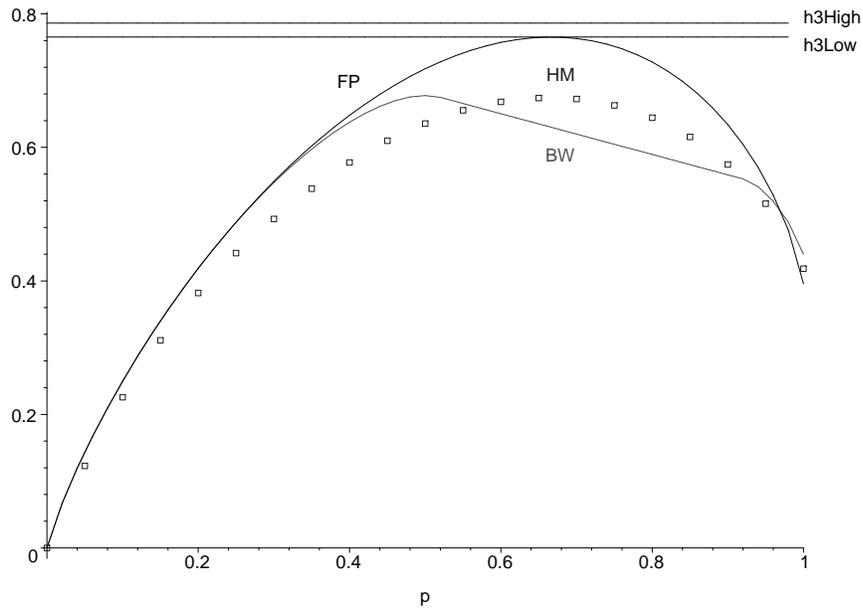}
  \caption{Lower bounds for $\lambda_3(p)$. HM is the lower bound of \cite{HaM}, BW is the lower bound of
  \cite{BoW}, FP is the lower bound of Theorem
  \protect\ref{lpmdczd}, h3Low and h3High are the best lower and
  upper bounds for $h_3 = \max \lambda_3(p)$.}\label{figbounds3}
 \end{figure}

 We conclude with a comparison of several lower bounds for the
 monomer-dimer entropy with dimer density $p$, $\lambda_d(p)$, for
 $d=2,3$. Hammersley and Mennon \cite{HaM} give a lower bound for
 $\lambda_d(p)$, graphed and tabulated in increments of $0.05$ for
 $ 0 \leq p \leq 1$. Bondy and Welsh \cite{BoW} give another lower
 bound for $\lambda_d(p)$, which depends on the dimer entropy
 $\lambda_d(1)$ and increases with it. Since $\lambda_3(1)$ is
 known only through upper and lower bounds, the bound of
 \cite{BoW} improves  each time a better lower bound for
 $\lambda_d(1)$ is found. We computed the lower bound of
 \cite{BoW} for $\lambda_3(p)$ using the best available lower
 bound $\lambda_3(1) = \widetilde{h}_3 \geq 0.4400758$. Hammersley
 and Mennon too tabulated and graphed the bound of \cite{BoW} for
 $\lambda_3(p)$, but at the time the available lower bound for
 $\lambda_3(1)$ was weaker. Figures \ref{figbounds2} and
 \ref{figbounds3} illustrate the lower bounds for $\lambda_d(p)$,
 $d = 2,3$, due to \cite{HaM}, \cite{BoW}, and Theorem
 \ref{lpmdczd}. Figure \ref{figbounds2} also illustrates the Monte
 Carlo estimates of \cite{HaM}. It is seen that except for very
 high $p$, the best lower bound is given by Theorem \ref{lpmdczd}.
 (As pointed out above, (\ref{cnvh2}) implies that the Monte Carlo
 estimates above the line $y=h_2$ are over estimates). We also
 include in the figure estimates of $\lambda_2(p)$ obtained from
 the heuristic computations of Baxter \cite{Bax1}. One can obtain
 from the lower bound of \cite{Wan} a corresponding lower bound
 for $\lambda_d(p)$. It turns out that for $d=2,3$, this bound is
 dominated by the maximum of the lower bound given by Theorem
 \ref{lpmdczd} and the lower bound of \cite{BoW}.

 \clearpage
 \bibliographystyle{plain}

\begin{thebibliography}{MMM}
 \bibitem{Bax1} R.J. Baxter, Dimers on a rectangular lattice, \emph{J.
 Math.\ Phys.} 9 (1968), 650--654.
 \bibitem{Bax}  R.J. Baxter, Hard hexagons: exact solution, \emph{J. Phys.\ A: Math.\ Gen.}
 13 (1980), L61--L70.
 \bibitem{BeS}  I. Beichl and F. Sullivan, Approximating the permanent via importance
 sampling with application to the dimer cover problem, \emph{J. Comp.\ Phys.} 149 (1999),
 128--147.
 \bibitem{Ber}  R. Berger, The undecidability of the domino problem, \emph{Mem.\ Amer.\ Math.\ Soc.}
 66 (1966).
 \bibitem{BoW}  J.A. Bondy and D.J.A. Welsh, A note on the monomer
 dimer problem, \emph{Proc.\ Camb.\ Phil.\ Soc.} 62 (1966), 503--505.
 \bibitem{CaW}  N.J. Calkin and H.S. Wilf, The number of independent sets in a grid graph,
 \emph{SIAM J. Discr.\ Math.} 11 (1998), 54--60.
  \bibitem{Ciu}  M. Ciucu, An improved upper bound for the 3-dimensional dimer problem, \emph{Duke Math.\ J.}
  94 (1998), 1--11.
 \bibitem{Fis}  M.E. Fisher, Statistical mechanics of dimers on a plane lattice, \emph{Phys.\ Rev.} 124
 (1961), 1664--1672.
 \bibitem{FoJ}  S. Forschhammer and J. Justesen, Entropy bounds for constrained two-dimensional random
 fileds, \emph{IEEE Trans.\ Info.\ Theory} 45 (1999), 118--127.
 \bibitem{FoR}  R.H. Fowler and G.S. Rushbrooke, Statistical theory of perfect solutions,
 \emph{Trans.\ Faraday Soc.} 33 (1937), 1272--1294.
 \bibitem{Fr3}  S. Friedland, A proof of a generalized van der Waerden conjecture on permanents,
 \emph{Lin.\ Multilin.\ Algebra} 11 (1982), 107--120.
 \bibitem{Fr1}  S. Friedland,
 On the entropy of Z-d subshifts of finite type, \emph{Linear Algebra Appl.} 252 (1997), 199--220.
 \bibitem{Fr2}  S. Friedland, Multi-dimensional capacity, pressure and Hausdorff dimension,
 in \emph{Mathematical System Theory in Biology,
 Communication, Computation and Finance}, edited by J. Rosenthal and D. Gilliam,
 IMA Vol.\ Ser.\ 134, Springer, New York, 2003, 183--222.
 \bibitem{HL} O.J. Heilman and E.H. Lieb, Theory of monomer-dimer
 systems, \emph{Comm.\ Math.\ Phys.} 25 (1972), 190--232; Errata
 27 (1972), 166.
 \bibitem{Gau}  D.S. Gaunt, Exact series-expansion study of the monomer-dimer problem, \emph{Phys.\
 Rev.} 179 (1969) 174--186.
 \bibitem{Ha1}  J.M. Hammersley, Existence theorems and Monte Carlo methods for the monomer-dimer
 problem, in \emph{Reseach papers in statistics: Festschrift for J. Neyman}, edited by F.N. David,
 Wiley, London, 1966, 125--146.
 \bibitem{Ha2}  J.M. Hammersley, An improved lower bound for the multidimesional dimer problem,
 \emph{Proc.\ Camb.\ Phil.\ Soc.} 64 (1966), 455--463.
 \bibitem{Ha3} J.M. Hammersley, Calculations of lattice statistics, in
 \emph{Proc.\ Comput.\ Physics Con.}, London: Inst.\ of Phys.\ \& Phys.\ Soc.,
 1970, 1--8.
 \bibitem{HaM}  J. Hammersley and V. Menon, A lower bound for the monomer-dimer problem,
 \emph{J. Inst.\ Math.\ Applic.} 6 (1970), 341--364.
 \bibitem{HKC} L.P. Hurd, J. Kari, and K. Culik, The
 topological entropy of cellular automata is uncomputable, \emph{Ergodic Theory Dynam.\ Systems}, 12 (1992), 255--265.
 \bibitem{Kas}  P.W. Kasteleyn, The statistics of dimers on a lattice, \emph{Physica} 27 (1961),
 1209--1225.
 \bibitem{KeRaSi} C. Kenyon, D.Randall and A. Sinclair,
 Approximating the number of monomer-dimer coverings of a lattice,
 J. Stat.\ Phys.\ 83 (1996), 637--659.
 \bibitem{Jer} M. Jerrum, Two-dimensional monomer-dimer systems are computationally intractible,
 \emph{J. Stat.\ Phys.} 48 (1987), 121--134.
 \bibitem{L}  E.H. Lieb, The solution of the dimer problem by the
 transfer matrix method, \emph{J. Math.\ Phys.} 8 (1967),
 2339--2341.
 \bibitem{Lie}  E.H. Lieb, Residual entropy of square ice, \emph{Phys.\ Review} 162 (1967), 162--172.
 \bibitem{Lun} P.H. Lundow, Compression of transfer matrices,
 \emph{Discrete Math.} 231 (2001), 321--329.
 \bibitem{Na2}  J.F. Nagle, New series expansion method for the dimer problem, \emph{Phys.\ Rev.} 152
 (1966), 190--197.
 \bibitem{NaZ}  Z. Nagy and K. Zeger, Capacity bounds for the 3-dimensional $(0,1)$ run length limited
 channel, \emph{IEEE Trans.\ Info.\ Theory} 46 (2000), 1030--1033.
 \bibitem{Pau}  L. Pauling, \emph{J. Amer.\ Chem.\ Soc.} 57 (1935), 2680--.
 \bibitem{Sc}  K. Schmidt, \emph{Algebraic Ideas in Ergodic Theory}, Amer.\ Math.\ Soc., 1990.
 \bibitem{Sch}  A. Schrijver, Counting $1$-factors in regular bipartite graphs,
 \emph{J. Comb.\ Theory} B 72 (1998), 122--135.
 \bibitem{Wan} I.M. Wanless, A lower bound on the maximum
 permanent in $\Lambda_n^k$, \emph{Linear Algebra Appl.} 373
 (2003), 153--167.
 \bibitem{WeB}  W. Weeks and R.E. Blahut, The capacity and coding gain of certain checkerboard codes,
 \emph{IEEE Trans.\ Info.\ Theory} 44 (1998), 1193--1203.


 \end{thebibliography}

 \end{document}